\documentclass[10pt,a4paper]{amsart}
\usepackage{amscd,amssymb,amsopn,amsmath,amsthm,graphics,amsfonts,enumerate,verbatim,calc}
\usepackage[dvips]{graphicx}

\usepackage{mathpazo}
\usepackage{color}
\usepackage{latexsym}
\usepackage{amsthm,amsfonts,amssymb,mathrsfs}
\usepackage{rotating}
\usepackage[leqno]{amsmath}
\usepackage{xspace}
\usepackage[all]{xy}
\usepackage{longtable}
\headsep=1cm \numberwithin{equation}{section}
\newtheorem{theorem}{Theorem}[section]
\newtheorem{lemma}[theorem]{Lemma}

\newtheorem{proposition}[theorem]{Proposition}
\newtheorem{corollary}[theorem]{Corollary}

\theoremstyle{definition}
\newtheorem{definition}[theorem]{Definition}
\theoremstyle{remark}
\newtheorem{remark}[theorem]{Remark}
\newtheorem{fact}[theorem]{Fact}
\newtheorem{example}[theorem]{Example}
\newtheorem{observation}[theorem]{Observation}
\newtheorem{discussion}[theorem]{Discussion}
\newtheorem{question}[theorem]{Question}
\newtheorem{conjecture}[theorem]{Conjecture}

\newtheorem{acknowledgement}{Acknowledgement}

\newcommand{\UFD}{\operatorname{UFD}}
\newcommand{\GCD}{\operatorname{GCD}}
\newcommand{\DVR}{\operatorname{DVR}}

\newcommand{\Ass}{\operatorname{Ass}}
\newcommand{\im}{\operatorname{im}}

\newcommand{\pgrade}{\operatorname{p.grade}}
\newcommand{\Kgrade}{\operatorname{K.grade}}
\newcommand{\ara}{\operatorname{ara}}
\newcommand{\cgrade}{\operatorname{c.grade}}
\newcommand{\Egrade}{\operatorname{E.grade}}
\newcommand{\Cgrade}{\operatorname{\check{C}.grade}}

\newcommand{\Spec}{\operatorname{Spec}}

\newcommand{\rad}{\operatorname{rad}}

\newcommand{\hens}{\operatorname{hens}}

\newcommand{\Ht}{\operatorname{ht}}
\newcommand{\pd}{\operatorname{p.dim}}
\newcommand{\fpd}{\operatorname{\mathfrak{f}in}}
\newcommand{\Fpd}{\operatorname{fin}}

\newcommand{\gd}{\operatorname{gl.dim}}

\newcommand{\F}{\operatorname{F}}
\newcommand{\e}{\operatorname{e}}
\newcommand{\nil}{\operatorname{nil}}
\newcommand{\rank}{\operatorname{rank}}

\newcommand{\Wdim}{\operatorname{w.dim}}
\newcommand{\fd}{\operatorname{fl.dim}}

\newcommand{\V}{\operatorname{V}}

\newcommand{\id}{\operatorname{id}}
\newcommand{\Ext}{\operatorname{Ext}}
\newcommand{\E}{\operatorname{\textbf{E}}}

\newcommand{\Tor}{\operatorname{Tor}}

\newcommand{\Hom}{\operatorname{Hom}}

\newcommand{\Ann}{\operatorname{Ann}}

\newcommand{\red}{\operatorname{red}}

\newcommand{\Jac}{\operatorname{Jac}}
\newcommand{\depth}{\operatorname{depth}}
\newcommand{\Kdepth}{\operatorname{Kdepth}}

\newcommand{\Char}{\operatorname{char}}
\newcommand{\coker}{\operatorname{coker}}
\newcommand{\HH}{\operatorname{H}}
\newcommand{\W}{\textbf{W}}
\newcommand{\trdeg}{\operatorname{tr.deg}}

\newcommand{\vpl}{\operatornamewithlimits{\varprojlim}}

\newcommand{\lo}{\longrightarrow}
\newcommand{\fm}{\frak{m}}
\newcommand{\fp}{\frak{p}}
\newcommand{\fq}{\frak{q}}
\newcommand{\fa}{\frak{a}}
\newcommand{\fb}{\frak{b}}
\newcommand{\fc}{\frak{c}}

\newcommand{\fn}{\frak{n}}

\begin{document}

\author[]{mohsen asgharzadeh}

\address{}
\email{mohsenasgharzadeh@gmail.com}

\title[ ]
{homological aspects of perfect algebras}


\begin{abstract}
We investigate homological properties of perfect algebras of prime characteristic. Our  principle is as follows: perfect algebras resolve  the singularities.
For example,  we show any module over the ring of absolute integral closure has finite flat dimension.
Under some mild conditions, we show any module over this ring has finite projective dimension.
We compute weak
and global dimensions of perfect rings in a series of nontrivial cases. Some interesting applications are given.
In particular, we answer some questions asked by Shimomoto.
\end{abstract}

\maketitle

\section{Introduction}

A commutative ring of   characteristic $p$ is called \textit{perfect} if the Frobenius map is an isomorphism. For many proposes, the surjectivity is enough.
Also,  a  ring $A$ of (mixed) characteristic $p$ is called  \textit{semiperfect} if its mod $p$ reduction  has all of its $p$-power roots.
Perfect  rings get interesting nowadays. Our interest in perfect algebras is as follows: perfect algebras resolve  the singularities.
 Let us recall some important examples of  perfect algebras. Let $R$ be a local domain. The ring of \textit{absolute integral closure} $R^\textbf{+}$ is
the integral closure of $R$ inside an algebraic closure
of the field of fractions of $R$. The symbol $R^\textbf{+}$ introduced by Artin in \cite{Ar}, where among other
things, he proved  $R^\textbf{+}$ has only one maximal ideal $\fm_{R^\textbf{+}}$,
when $(R,\fm)$ is   henselian. By using an  idea due to  Bhatt and Scholze  \cite{BS} we  show:

\begin{theorem}  \label{fill}
Let $(R,\fm,k)$ be a
local complete domain of  prime characteristic. Then any $R^\textbf{+}$-module  has finite flat dimension.
In addition, if $k$ is of $\aleph_n$-cardinality (e.g., $k$ is countable), then any $R^\textbf{+}$-module has finite projective dimension.
\end{theorem}

 The first part extends a result of Hochster and Aberbach \cite{ab} to the full setting.
Computing homological dimensions over $R^\textbf{+}$  is  difficult when $R$ is of mixed characteristic (even for simple modules). Here is a sample:

 \begin{fact}\label{ab}(See \cite[Theorem 3.5]{ab})
$\fd_{R^\textbf{+}}(R^\textbf{+}/ \fm_{R^\textbf{+}})=\dim (R)$ in all mixed characteristic cases if and only
 the direct summand conjecture holds in mixed characteristic.
 \end{fact}

The direct summand conjecture is now a theorem by Andr\'{e}'s recent work \cite{an}.
The organization of the paper is as follows. In
\S 2, we give a quick review of homological invariants such as
weak dimension, global dimension, and depth. In the sequel we will use all of them. A
landmark  result on these invariants is due to Auslander and Serre:  if a ring $(A,\fm)$ is noetherian and local then $$\depth(A)\leq\dim (A)=\ara(\fm)\leq\mu(\fm)\leq\Wdim (A)=\gd(A).$$Morever, the equality holds if $\gd(A)$ is finite (here, $\ara(\fm)$ is the minimum number of elements required to generate
$\fm$ up to radical).
The situation is not so simple if $A$ is not noetherian.
For example, the ring of entire functions  has infinite Krull dimension
and finite global dimension (see \cite{O2}).
We will use the following result to show certain perfect algebras are not coherent:

\begin{observation}\label{wdim=dime}
Let $(A,\fm)$ be a coherent quasilocal  ring and  of finite weak dimension. Then $\Kdepth (A)=\Wdim (A) \leq\dim (A).$ In particular, $\Wdim (A)\leq\ara(\fm)$.
\end{observation}

As an application, we reprove a theorem of Vasconcelos (see Corollary \ref{wdim=emb}).

 \S 3 is devoted to computing  homological invariants over $R^\infty$. Following Greenberg \cite{g} and Serre \cite{Ser}, the symbol $R^\infty$ stands
for the perfect closure of  a noetherian ring $R$. As far as we know, Vasconcelos (was the first person who) computed the weak dimension of a very special ring of the form $R^\infty$, see \cite[Example 5.28]{V2}. Revisiting \cite{BS}, we observe that  $\gd(R^{\infty})$ is finite. In fact the following stated in \cite[Footnote 24]{BS} without a proof.

\begin{corollary}\label{g}
Let $(R,\fm, k)$ be a complete local  domain  of prime characteristic. Then $\gd(R^\infty)\leq2\dim( R)+1$.
\end{corollary}

This corollary answers \cite[Question 9.2]{A} and \cite[Question 9.5(ii)]{A}.
The structure of free resolutions over $R^\infty$ is quite mysterious. However
 for radical ideals, we compute some explicit free resolutions:

\begin{observation}(See Proposition \ref{str})
For any radical ideal  $\fa$ of $R^{\infty}$,  the module $R^{\infty}/ \fa$ has a   free resolution of countably generated
free $R^{\infty}$-modules of length  bounded by $2\dim R$.
 \end{observation}

This extends \cite[Theorem 1.2(i)]{A} by presenting the bound $2\dim R$.
The next problem is as follows: What is $\gd(R^{\infty})$? We answer this in a low-dimensional case, see  Corollary \ref{1}.
Over complete regular local rings that
are not  field,  we show $$\dim (R^\infty)+1=\gd(R^{\infty})>\Wdim (R^{\infty})=\Kdepth (R^{\infty})=\dim (R^\infty)=\ara(\fm_{R^{\infty}}).$$
Recall that $R$ is called \textit{$F$-coherent} if $R^\infty$ is coherent. In \S 4 we deal with the following questions asked by Shimomoto:

\begin{question}\label{introsh}Let $R$ be a local ring of prime characteristic and $t\in R$  a non-zero divisor.
\begin{enumerate}
\item[i)] (See \cite[Question 2]{ss})   Let $R/t R$ be $F$-coherent. Is $R$ $F$-coherent?
\item[ii)] (See \cite[Question 3]{ss})  Let $R$ be $F$-coherent. Is the Hilbert--Kunz multiplicity of $R$ rational?
\item[iii)] (See \cite[Question 1]{ss})  Let $R$ be $F$-coherent and $R\to S$ be flat. What conditions on the fibers required to $S$ be
$F$-coherent?
\end{enumerate}
\end{question}
We done Question \ref{introsh}(i) by presenting a perfect ring with zero-divisors such that $$\dim (A)=\Kdepth (A)=1<\Wdim (A)=2<\gd(A)=3.$$In particular, there is ``a commutative local ring with finite global dimension and zero divisors.'' Also, see   \cite{O2}. Remark \ref{newex}
resolves Question \ref{introsh}(i) by another method.
Concerning Question \ref{introsh}(ii), we ascend up to the perfection and descend down to $R$. I.e., we have:

\begin{proposition}
Let $(R,\fm,k)$ be an $\F$-finite and  $\F$-coherent domain  and $I\lhd R$ be  of finite colength. If $R$ is  Cohen-Macaulay, then $\e_{HK}(I,R)$ is rational.
\end{proposition}

Suppose $R\to S$ is of finite presentation. Then
 $R\to S$ is  \'{e}tale if and only if it is flat and unramified. On the other hand unramified property defined only by the study of fibers.
 One can drop the finiteness by looking at  \textit{weakly \'{e}tale} extensions. In particular, the following partially answers Question \ref{introsh}(iii):
Let $A$ be a  local $\F$-coherent and $B$ be weakly \'{e}tale over $A$. Then $B$ is $\F$-coherent.

  In \S 5, we use the results of \S 3 to prove Theorem  \ref{fill}.
\S 6, deals with \textit{desingularization}  of $R^\textbf{+}$.
As a corollary  to the presented results, we show:

\begin{corollary}
Let $R$ be a complete local domain of mixed characteristic. If
$\dim (R )> 3$ then $R^\textbf{+}$ is not coherent.
\end{corollary}
 This is the  mixed characteristic
version of a result of Hochster and Aberbach.
In \S 7, we present situations for which the global dimension of certain perfect
algebras depend on the characteristic:
\begin{example}
Let $R:=\mathbb{F}_p[t,t\sqrt{t+1}]$. Then
\begin{equation*}
\gd(R^{\infty})= \left\{
\begin{array}{rl}
3 & \  \   \   \   \   \ \  \   \   \   \   \ \text{if } p\neq 2\\
2 & \  \   \   \   \   \ \  \   \   \   \   \ \text{if } p=2
\end{array} \right.\  \ \ \ \ \ \ \ \
\Wdim(R^{\infty})= \left\{
\begin{array}{rl}
2 & \  \   \   \   \   \ \  \   \   \   \   \ \text{if } p\neq 2\\
1 & \  \   \   \   \   \ \  \   \   \   \   \ \text{if } p=2
\end{array} \right.
\end{equation*}
\end{example}

In \S 8, we present a simple proof of a \textit{miraculous vanishing formula}  due to Bhatt and Scholze. We drive it
from the special case presented in our pervious work \cite{A}, where it related to the so called telescope conjecture. This has some applications.
We finish \S 8 by the following result:

\begin{observation}
Let $R \to S$ be a perfectly finitely presented map of perfect $\mathbb{F}_p$-algebras. Then $\pd_R(S)<\infty$.
\end{observation}

Finiteness of  $\fd_R(S)<\infty$ is a subject of a   recent result of  Bhatt and Scholze.
In \S 9 we answer a question asked by Shimomoto (see \cite[Question 2]{sss}):

\begin{example}(After Kedlaya)
Let $R$ be any $1$-dimensional  local ring of prime characteristic. Then $\W(R^\infty)$ is not coherent
without any regard with respect to coherent property of $R^\infty$.
\end{example}

Also, the following   extends \cite[Lemma 7.8]{BS} by Bhatt and Scholze:

\begin{remark}
 Let $R$ be perfect  and $Q$ be (not necessarily finitely generated and not necessarily projective) an  $R$-module of finite projective dimension.
Then $\pd_{\W(R)} (Q)=\pd_{R} (Q)+1.$
\end{remark}

We emphasize that
perfect algebras are almost non-noetherian. This is the main difficulty.
Despite of this, perfect closure of a noetherian ring with  singularity is a ring of finite global dimension.

\section{homological invariants}

In this note all rings are  commutative.  We consider  both
noetherian and  non-noetherian rings.
Recall that by $\pd(-)$ (resp. $\fd(-)$), we mean the
projective dimension (resp. the flat dimension).
A quasilocal ring  $A$ is a commutative ring with a unique maximal ideal $\fm_A$. A  local
ring  is a noetherian quasilocal ring.
By $\gd(-)$  we mean the \textit{global dimension}.
Also, $\Wdim(-)$ stands for the \textit{weak dimension}.  Recall  for any commutative ring $A$ that
\[\begin{array}{ll}
\Wdim(A)&:=\sup\{\fd(M):M \emph{ is an $A$-module }\}\\&=\sup\{\fd(A/
\fa):\fa \emph{ is a finitely generated ideal of }A\}\quad(2.1.1)
\\
\end{array}\] In particular, if flat dimension of any finitely generated ideal
is bounded by a uniform integer $n$, then flat dimension any module is bounded by the integer $n$.

Recall that a set $\Gamma$ is an
ordinal if $\Gamma$ is totally ordered with respect to inclusion and every
element of $\Gamma$ is  a subset of $\Gamma$. Also, one can see that
$\Omega$ is itself an ordinal number larger than all countable ones,
so it is an uncountable set. By $\aleph_{-1}$ we denote the cardinality of
finite sets. By $\aleph_0$ we mean the cardinality of the set of all
natural numbers. We look at $$\Omega:=\{\alpha:\alpha
\textit{ is a countable ordinal number} \}.$$   By definition
$\aleph_1$, is the cardinality of $\Omega$. Inductively, $\aleph_n$
can be defined for all $n\in\mathbb{N}$. A ring is
called $\aleph_n$-noetherian if each of its ideals  can be generated
by a set of cardinality bounded by $\aleph_n$. So, noetherian rings
are exactly $\aleph_{-1}$-noetherian rings.

\begin{lemma}\label{j} (See the proof of \cite[Corollary 2.47]{O1})
Let  $\fa$ be an ideal of a $\aleph_n$-noetherian ring
$A$. Then $\pd_{A}(A/\fa)\leq\fd_{A}(A/\fa)+n+1$.
\end{lemma}

A ring is called \textit{coherent}, if its finitely generated ideals are finitely presented.  The following
is a way to show a ring is coherent.

\begin{fact}\label{flatdir}
(See \cite[Theorem 2.3.2]{G}) Any flat direct limit of coherent rings is coherent.
\end{fact}

\begin{definition}\label{kd}
The \textit{Koszul grade} of a finitely
generated ideal  with a generating set $\underline{x}$ on  a module $M$
is defined by
$$\Kgrade_A(\underline{x},M):=\inf\{i \in\mathbb{N}\cup\{0\} | H^{i}(\Hom_A(
K_{\bullet}(\underline{x}), M)) \neq0\}.$$For
an ideal $\frak b$ (not necessarily finitely generated), \textit{Koszul grade} of
$\fa$ on $M$ can be defined by
$$\Kgrade_A(\fb,M):=\sup\{\Kgrade_A(\frak c,M):\frak c\in\Sigma\}\quad(\ref{kd}.1)$$ where $\Sigma$ is the family of all finitely generated
subideals  of $\fa$.   The notation  $\Kdepth(-)$ stands for $\Kgrade_A(\fm,-)$ where  $(A,\fm)$ is quasilocal.
\end{definition}

\begin{remark} \label{p=k}
i) (See \cite[Page 149]{No}) The classical grade  of $\fa$ on $M$,
denoted by $\cgrade_A(\fa,M)$, is defined to the supremum of the
lengths of all weak regular sequences on $M$ contained in $\fa$. The polynomial grade of $\fa$ on $M$ is
defined by
$$\pgrade_A(\fa,M):=\underset{m\rightarrow\infty}{\lim}
\cgrade_{A[t_1, \cdots,t_m]}(\fa A[t_1, \cdots,t_m],A[t_1,,
\cdots,t_m]\otimes_A M).$$

ii) One has $\pgrade_A(\fa,-)=\Kgrade_A(\fa,-)$, see e.g. \cite[Proposition 2.3]{AT}.
\end{remark}

Let us cite the following basic properties of Koszul grade.

\begin{fact}\label{cp} Let $R$ be any ring, $M$ an $R$-module and $\fa$ an ideal. The following holds:\begin{enumerate}
\item[i)] One has $\Kgrade_{R}(\fa,M)=\Kgrade_{R}(\fp,M)$ for some prime ideal $\fp$ (see \cite[Theorem 5.16]{No}).
\item[ii)] $\Kgrade$ is unique up to radical by \cite[Proposition 2.2 (vi)]{BH}.
\item[iii)] (See \cite[Proposition 9.1.2(g)]{BH}) If $S\subset R$ containing a system of generators $\underline{x}$ of $\fa$ then
$$\Kgrade_{R}(\fa,M)=\Kgrade_{S}(\underline{x}S,M).$$
\item[iv)] (See \cite[Proposition 9.1.4]{BH}) $\Kgrade(\fa,R)=\inf\{\Kdepth_{R_{\fp}}(R_{\fp}):\fp\in V(\fa)\}$.
\item[v)] (See \cite[Theorem  9.1.6]{BH}; Buchsbaum-Eisenbud, Northcott)
Let $$\textbf{F}:\xymatrix{0\ar[r]&F_{m}\ar[r]&\ldots\ar[r]&F_{j+1}\ar[r]^{f_j}&F_{j}\ar[r]&\ldots\ar[r]&F_{0}\ar[r]&0,}$$
be a complex of finite free $R$-modules
and $r_i$ be the expected rank of $f_i$. Then
$\textbf{F}\otimes_RM$ is acyclic
 if and only if  $\Kgrade_R(I_{r_i}(f_i),M)\geq	 i$ for all $i$.
\item[vi)](Auslander-Buchsbaum, Hochster \cite[Chap. 6, Theorem 2]{No})
Suppose $\textbf{F}$ in the above item is acyclic and $R $ is quasilocal. Let $N:=\coker(f_0)$. Then
$$\Kdepth_R(N)+\pd_R(N)=\Kdepth_R(R).$$
\end{enumerate}
\end{fact}

\begin{theorem}\label{wdim=dim}
Let $(A,\fm)$ be a coherent quasilocal  ring of finite weak dimension. Then $\Wdim (A)=\Kdepth (A) \leq\dim (A).$ In particular, $\Wdim (A)\leq\ara(\fm)$.
\end{theorem}

\begin{proof}
Let $\fa\vartriangleleft A$ be finitely generated. Then $A/\fa$ is finitely presented and is of finite flat dimension.
Finitely present flat modules over quasilocal rings are  free. It turns out that $$\fd(A/ \fa)=\pd(A/ \fa) \quad(\dagger).$$ Thus,
$A/\fa$  has finite free resolution by finitely generated free modules (see \cite[Corollary 2.5.2]{G}).
By Auslander-Buchsbaum-Hochster Fact \ref{cp}(vi), \[\begin{array}{ll}
\fd(A/ \fa)&=\pd(A/ \fa)\\
&=\Kdepth(A)-\Kdepth(A/ \fa)\\
&\leq\Kdepth(A).
\end{array}\]From this we deduce that $$\Wdim (A)\leq\Kdepth (A) \quad(\ddagger)$$
\begin{enumerate}
\item[Fact A:](See \cite[Lemma 3.2]{AT}) Let $\fa$ be an ideal of a ring $B$ and $M$ a finitely generated
$B$-module. Then $$\Kgrade_B (\fa,M)\leq \Ht_{M}(\fa).$$
\end{enumerate}

In view of Fact A, $\Kdepth(A)\leq\dim (A)$. Thus, $\Wdim (A)\leq \dim (A) $.
To show $\Wdim (A)=\Kdepth (A)$ we need to recall the concept of Ext-grade. The $\Ext$ grade of $\fa$ on $-$ is defined by
$$\Egrade_{A}(\fa,-):=\inf\{i\in \mathbb{N}\cup\{0\}|\Ext^{i}_{A}(A/\fa,
-)\neq0\}.$$In general $\Egrade_{A}(\fa,-)\neq\Kgrade_{R}(\fa,-)$\footnote{let
$R:=\mathbb{Q}[x_n: n\in \mathbb{N}]/(x_n^n: n\in \mathbb{N})$. Set $\fa:=(x_n: n\in \mathbb{N})$.
By \cite[Page 367]{B},    $\Kgrade_R(\fa,R)\neq
\Egrade_R(\fa,R).$ }. However, if $\fa$ is finitely generated $$\Egrade_{A}(\fa,-)=\Kgrade_{A}(\fa,-) \quad(\natural)$$ (see \cite[Proposition 2.3(iii)]{AT}).  Clearly $$\Egrade_{A}(\fa,-)\leq \pd(A/\fa)\quad(\diamond)$$ Let $\Sigma$ be the family of all finitely generated
subideals  of $\fm$. Thus, \[\begin{array}{ll}
\Wdim (A)&\stackrel{(\ddagger)}\leq\Kgrade_A(\fm,A)\\
&\stackrel{\ref{kd}.1}=\sup\{\Kgrade_A(\frak a,A):\frak a\in\Sigma\}\\
&\stackrel{(\natural)}=\sup\{\Egrade_A(\frak a,A):\frak a\in\Sigma\}\\
&\stackrel{(\diamond)}\leq\sup\{\pd(A/\fa):\frak a\in\Sigma\}\\
&\stackrel{(\dagger)}=\sup\{\fd(A/\fa):\frak a\in\Sigma\}\\
&\stackrel{2.1.1}=\Wdim (A).
\end{array}\] Thus $\Wdim (A)=\Kdepth(A)$. In order to show
 $\Wdim (A)\leq\ara(\fm)$ we may assume that $\ell:=\ara(\fm)<\infty$. Let $\underline{x}:=x_1,\ldots,x_{\ell}$ be such that $\rad(\underline{x})=\fm$. In view of Fact \ref{cp} $\Kgrade(\fm,A)=\Kgrade(\underline{x},A)$.
By definition, $\Kdepth(A)=\Kgrade(\fm,A)=\Kgrade(\underline{x},A)\leq \ell$. By the first part, $\Wdim (A)=\Kdepth(A)\leq\ara(\fm).$
\end{proof}

In the  proof of Theorem \ref{wdim=dim} the following invariant appeared:
by very small finitistic dimension we mean
$$\fpd(A):=\sup\{\pd(A/\fa):\fa\textit{ is finitely generated and of finite projective dimension}\}.$$
 Recall that the classical small finitistic dimension is
$$\Fpd(A):=\sup\{\pd(M):M\textit{ is finitely generated and of finite projective dimension}\}.$$
It is easy to find examples with $\fpd(A)\lvertneqq\Fpd(A)$: Let $A:=(\mathbb{F}_2[[X]])^\infty$. Then   $\fpd(A)=1\lvertneqq\Fpd(A)=2$.
Also there is a situation for which $\Wdim(A)\lvertneqq\fpd(A)$ (such a ring is not coherent):

\begin{example}
 Let $A$ be the subring of $C(\mathbb{R})$ (the ring of all
continuous real-valued functions) consisting of piecewise
sums of odd roots of polynomials and quotients thereof. By \cite{O2} $\Wdim(A)=2$.  Also, Osofsky \cite{O2}
presents an element $f\in A$ such that $\pd(A/fA)=3$. So, $\Wdim(A)\lvertneqq\fpd(A)$.
\end{example}

\begin{lemma}\label{should1}
 Let $(A,\fm)$ be a coherent  quasilocal ring  and $\underline{x}:=x_1,\ldots,x_d\subset\fm$ be such that $\Kgrade (\underline{x},A)=d$. Then $\underline{x}$ is a regular sequence.
\end{lemma}

\begin{proof}
Let $1\leq i< d$ and set $\fa_i:=(x_1,\ldots,x_i)$. Since $A$ is coherent,
$H^{j}(\Hom_A( K_{\bullet}(\fa_i, A)))$ is
finitely generated.   By using  Nakayama's Lemma and an easy
induction  we see  that $\Kgrade_A(\fa_i,A)=
i$. In particular, $\Kgrade_A(x_1,A)=
1$. Thus, $x_1$ is a regular sequence. We note   that $R/x_1R$ is  coherent (see \cite[Theorem 2.4.1(1)]{G}). By using an easy induction on $d$
we deduce  that $\underline{x}$ is a regular sequence.
\end{proof}

\begin{corollary}\label{should}
 Let $(A,\fm)$ be a quasilocal coherent  ring such that $\fm$ is a radical of a finitely generated ideal with generating set $\underline{x}:=x_1,\ldots,x_d$ and
$\Kdepth (A)=d$. Then any permutation of $\underline{x}$ is a regular sequence over $A$.
\end{corollary}

\begin{proof}
By Fact \ref{cp} $\Kgrade_A(\underline{x},A)=\Kdepth(A)=d$.
In view of  Lemma \ref{should1}, $\underline{x}$ is a regular sequence. Since Koszul homology is invariant under permutation,
any permutation of $\underline{x}$ is a regular sequence.
\end{proof}

By $\mu(\fm)$  we mean the minimal number of
elements of $A$ that need to generate $\fm$. Here we reprove (and extend)
a result of Vasconcelos by a different argument (see \cite[Theorem 5.22]{V2}):

\begin{corollary}\label{wdim=emb}(Northcott+Vasconcelos)
Let $(A,\fm)$ be a  quasilocal  ring  and $\fm$ is finitely generated. Then $\mu(\fm)\leq\Wdim(A)$.  Suppose in addition that $A$ is coherent and  $\Wdim(A)<\infty$. Then $\fm$ is generated  by a regular sequence $\underline{x}$ of length $\Wdim(A)$. Any permutation of $\underline{x}$ is a regular sequence.
\end{corollary}

To find maximal regular sequences of different length  see \cite[Remark 5.23]{V2}.

\begin{proof}
In the light of  \cite[Theorem 3]{N} we see $\rank_{A/ \fm}(\frac{\fm}{\fm^2})\leq\Wdim(A).$ Since $\fm$ is finitely generated and by Nakayama's lemma, $d:=\mu(\fm)\leq \Wdim(A)$.  By definition, $\Kdepth(A)=\Kgrade(\fm,A)\leq \mu(\fm)$. Suppose  that $A$ is coherent and  of finite weak dimension. By  Theorem \ref{wdim=dim}, $$\mu(\fm)\leq \Wdim(A)=\Kdepth(A)\leq\mu(\fm).$$ Let $\underline{x}:=x_1,\ldots,x_d$ be a generating set of $\fm$. Due to Corollary \ref{should},  any permutation of $\underline{x}$ is a regular sequence.
\end{proof}

Recall from \cite{Ber} that a ring is  \textit{regular} if each   finitely generated ideal has finite projective dimension.
A  coherent  quasilocal ring  is called \textit{super regular}  if  its global dimension
is finite and equal to its weak dimension.
The following is due to Vasconcelos
and plays a  role in this paper:

\begin{fact}\label{sr}(See \cite[Theorem 5.29]{V2})
Let $(R,\fm)$ be a  super regular ring. Then $\fm$ can be generated by a regular sequence.
In particular, $\fm$ is finitely generated.
\end{fact}

 The extension $A \to B$ is called weakly \'{e}tale (or \textit{absolutely flat}) if $A \to B$ and $B\otimes_AB\to B$ are  flat.
The following  result is due to Olivier:
\begin{fact}\label{ol}(See \cite[Corollary 1]{Ol})
Let $A \to B$ be weakly \'{e}tale. Then $\Wdim(B)\leq\Wdim(A)$.
\end{fact}

\section{homological dimension over $R^\infty$ }

Rings in this section all are of prime characteristic $p$. Let $\F:R\to R$ be the Frobenius map. This sends $x$ to $x^p$. As an easy (but extremely important) fact, $\F$ is a ring homomorphism.
\begin{definition}
A ring of prime characteristic  is called perfect if the Frobenius map is an isomorphism.
\end{definition}

\begin{remark}
For many proposes the surjectivity of the Frobenius map is enough. Let us call such a ring as a \textit{semi-perfect} ring. Semi-perfect does not imply the perfectness. Note that
if  a ring is noetherian then any surjective ring-homomorphism is injective (see \cite[Ex. 3.6]{Mat}). But the noetherian assumption is important.
Because there are rings such as $A$ such that the Frobenius map over them is  surjective but not injective. The point is that
semi-perfect rings are (almost) non-noetherian (see the following observation for the explicit examples).
\end{remark}

By $F(R)$, we mean $R$ as a group equipped with  left and right scalar multiplication from $R$
given by
$a.r\star b = ab^pr,$ where $a,b\in R$ and $r\in F(R)$. Also, $\F^n(-):=(-)\otimes_R\F^n(R)$ is the \textit{Peskine-Szpiro} functor, please see \cite{PS2}.

\begin{definition}(Serre- Greenberg) Recall from \cite{G} that the  perfect closure $R^\infty$  of  $R$ is defined by $$R^\infty:=\varinjlim\xymatrix{(R\ar[r]^{\F}&R \ar[r]^{\F}&\ldots)}.$$
\end{definition}
Since $\F^{\gg 0}$
kills nilpotent elements, $R^\infty$ is reduced and exists uniquely. This sometimes is called the minimal perfect algebra.
Set $R_{\red}:=\frac{R}{\rad(0)}$. In fact, $R^\infty$ is defined
by adjoining to $R_{\red}$ all $p$-power roots of elements of $R_{\red}$.

\begin{fact}(Greenberg) Let $f:R\to S$ be a ring of prime characteristic $p$, there is a ring homomorphism $f^\infty:R^\infty\to S^\infty$
If $x\in R^\infty$, then $x^{p^n}\in R$ for some $n$. The assignment $x\mapsto f(x^{p^n})^{1/{p^n}}$ defines the well-defined map $f^\infty$. This makes the perfection as a functor.
In fact, Greenberg defined perfection of schemes.
\end{fact}

Recall from \cite{ss} that a ring is called $F$-coherent if its perfect closure is coherent, and we call $R^\infty$
 as  a coherent perfect closure  of $R$. Let us collect some elementary properties of perfect algebras that we need.

\begin{observation} \label{ob3} Let $(R,\fm)$ be a quasilocal ring of prime characteristic $p$. Then
\begin{enumerate}
\item[i)] Suppose $R$ is  noetherian. Then $R^\infty$ is noetherian if and only if $\dim R=0$.
\item[ii)] If $R$ is coherent and regular, then $R^\infty$ is coherent and flat over $R$.
\item[iii)] The coherent assumption in part ii) is really needed.
\item[iv)] The class of  $F$-coherent rings is strictly larger than the class of noetherian regular  rings.
\item[v)] Product of  radical ideals in a perfect algebra  is  the intersection of them.
\item[vi)] perfect algebras are  seminormal.
\item[vii)] Perfect algebras are not necessarily normal.
\item[viii)] Coherent perfect closure of a complete local domain is normal.
\item[ix)] A perfect domain  is $\UFD$ if and only if it is a field.
\item[x)] Any tensor product of perfect algebras is semi-perfect (e.g. has $p$-power root) but not necessarily perfect (e.g. roots are not unique).
Any tensor product of perfect algebras over a perfect ring is perfect.
For example,
$$(R_1\otimes_{R_0}R_2)^\infty\simeq R_1^\infty\otimes_{R_0^\infty}R_2^\infty.$$ Also, $\mathbb{F}_p$-endomorphism ring of a perfect algebra  has  $p$-power roots (possibly noncommutative).
\item[xi)] Any localization, direct limits, inverse limits and adic-completion of perfect algebras is again perfect. For example,
$(R^\infty)_{\fp}\simeq (R_{\fp\cap R})^\infty$ for any $\fp\in\Spec(R^\infty)$.
\item[xii)] If $A \to B$ is  weakly \'{e}tale, then $B^{1/p^n}\simeq B\otimes_AA^{1/p^n}$. In particular,
$B^\infty\simeq B\otimes_AA^\infty$.
\item[xiii)] Quotient of a perfect ring  by a radical ideal $\frak r$ is perfect. For example, $\frac{R^\infty}{\fp}\simeq (\frac{R}{\fp\cap R})^\infty$ for any $\fp\in\Spec(R^\infty)$.
\end{enumerate}
\end{observation}

\begin{proof}
\begin{enumerate}
\item[i)] If $R$  is zero dimension, then $R_{\red}$ is a field and its perfect closure is again a field. If $\dim R>0$, we take $x\in \fm$ which is not nil  and look at the increasing sequence $$0\subsetneqq (x)\subsetneqq (x^{1/p})\subsetneqq (x^{1/p^2})\subsetneqq\ldots \quad(\times)$$

\item[ii)] By a famous result of Kunz \cite{ku} (in the noetherian case), $R^{1/p}$ is flat over $R$. Let us show this in the coherent case as an application of the notion of Koszul grade.   We show  that $\Tor^R_i(R/\fa, F(R))=0$ for all $i > 0$ and for all finitely generated ideals $\fa\subset R$. Note that $R/\fa$  has a  free resolution $(F_\bullet,d_\bullet)$ consisting of finitely generated modules, since $R$ is coherent.  Then $(F_\bullet,d_\bullet)\otimes_R F(R)=(F_\bullet,d_\bullet^p)$. By $I_{t}(a_{ij})$ we mean the ideal
generated by the $t\times t$ minors of a matrix $(a_{ij})$.   Let $r_i$ be the expected rank of $d_\bullet$. Clearly,  $r_i$ is
the expected rank of $d_\bullet^p$. By Fact \ref{cp}(ii) $$\Kgrade_R(I_{r_i}(d_i), R)=\Kgrade_R(I_{r_i}(d_i^p ),R).$$
We apply  Fact \ref{cp}(v) (two times) to deduce that $(F_\bullet,d_\bullet^p)$ is exact. Hence $\Tor^R_i(R/\fa, F(R))=0$. So, $R^{1/p}$ is $R$-flat.
Therefore, $R^\infty$ is a flat directed union of coherent rings.
In view of Fact \ref{flatdir} $R^\infty$ is coherent.

\item[iii)] Let $A$ be a non $F$-coherent ring
(such a thing exists, see e.g. Example  \ref{123} below). We will see in Corollary \ref{gab} that $R:=A^\infty$ is regular.
By definition, $R^\infty=(A^\infty)^\infty=A^\infty$  which is not coherent.

\item[iv)] For example, the ring $\mathbb{F}_2[[X^2,XY,Y^2]]$ is $F$-coherent but not regular.
\item[v)] Prove this as an easy (but important) exercise or look at \cite[Proposition 2.11]{H1}.
\item[vi)]  Such a ring is reduced. In this case seminormality means that if $x\in Q(A)$ (the total quotient ring of $A$) is such that $x^2$ and $x^3$ are in $A$ then
$x\in A$. If $p=2$ or $p=3$ there is nothing to prove. Let us assume $p>3$. Then $p-3\in2\mathbb{N}$, e.g., $p=2r+3$ for some $r\in \mathbb{N}$. Thus,
$x^p\in A$. Since $A$ is perfect, we have $x\in A$ as claim.
\item[vii)] Any ring such that its normalization is not purely inseparable extension works. We left the details to the reader.
\item[viii)]  This is due to Shimomoto (see \cite[Theorem 3.8]{ss}).
\item[ix)] Remark that any $\UFD$ satisfies in the ascending chain condition on principal ideals. In view of $(\times)$ we get  the claim.
\item[x)]  Let $A_1$ and $A_2$ be perfect. Let $x:=\sum_{i=1}^m x_i\otimes y_i\in A_1\otimes_R A_2$. Set $
y:=\sum_{i=1}^m x_i^{1/p}\otimes y_i^{1/p}\in A_1\otimes_R A_2$. Then $y^p=x$. So, $A_1\otimes_R A_2$ has $p$-power root. Note that $R$ is not necessarily perfect. But the root is not necessarily unique.
For example, $a:=1\otimes x^{1/p}-x^{1/p}\otimes1\in (\mathbb{F}_p[[x]])^\infty\otimes_{\mathbb{F}_p[[x]]}(\mathbb{F}_p[[x]])^\infty$ is nonzero but is $p$-power is zero.
So, the Frobenius is not injective. 

Now, we assume that  $R$ is  perfect (i.e. the $p$-roots is unique). The claim follows by the following more general fact:
(If $Z\leftarrow X\rightarrow Y$ are perfectiod spaces over  a non-Archimedean field,
then $Z\times_XY$ exists in the category of perfectoid adic spaces. This interesting result is due to Scholze \cite[Proposition 6. 18]{Sch}.)
This implies that  $$(A_1\otimes_RA_2)^\infty=A_1^\infty\otimes_{R^\infty}A_2^\infty.$$ Thus, if all of $R$
and $A_1, A_2$ are perfect then all $p$-power roots of $A_1\otimes_RA_2$ is unique. So, $A_1\otimes_RA_2$ is perfect.

Let $A:=\Hom_{\mathbb{F}_p} (A_1,A_1)$. This is an associative ring. Let $f\in A$. Define $f^{1/p}(a):=(f(a))^{1/p}$.
It is easy to see it is additive. Let $r\in \mathbb{F}_p$ and $a\in A$. By Fermat's little theorem, $r=r^{p}$.\footnote{Associative rings that satisfy the polynomial identity $x^{n(x)}=x$ for some $n_x>1$, are very special: they are commutative. This is an interesting result of Jacobson.} Taking $p$-th root, we have $r^{1/p}=r$. Then $$f^{1/p}(ra)=(f(ra))^{1/p}=(rf(a))^{1/p}=r^{1/p}f(a)^{1/p}=rf(a)^{1/p}.$$
Thus, $f^{1/p}$ belongs to $A$ and that $(f^{1/p})^p=f$. So, $A$ has $p$-power root.
\item[xi)] This is  easy, and we left it to the reader.
\item[xii)] This is in \cite[Theorem 3.5.13]{almost} by Gabber and Ramero.
\item[xiii)] Let $r+\frak r\in A/\frak r$. Define $(r+\frak r)^{1/p}:=r^{1/p}+\frak r$. This is well-defined, because $\frak r$ is a radical ideal. The particular claim is also trivial.
\end{enumerate}
\end{proof}

\begin{remark}
Recall that a $\GCD$ domain is an integral domain  with the property that any two non-zero elements have a greatest common divisor. This is well-known  that
any $\GCD$  domain is normal.
\end{remark}

Homological properties of perfect algebras not only simplified things but also extend them:

\begin{corollary} \label{ob2}  (Compare with Observation \ref{ob3}(viii)) Coherent perfect closure of a local domain $R$ is a $\GCD$ domain. In particular, $R^\infty$ is
 a normal domain.
\end{corollary}

\begin{proof}
Recall that $R^\infty$ is quasilocal, coherent and regular (see \cite[Theorem 1.2(iii)]{A}). By \cite[Corollary 6.2.10]{G} any quasilocal, coherent and regular domain is a $\GCD$ domain. By the above remark, $R^\infty$  is normal.
\end{proof}

\begin{fact}(See \cite[Proposition 11.31]{BS}) \label{uni}
Let $R^\infty$ be the perfection of a complete
local ring $R$. Then $R^\infty$ has finite global dimension. In fact,  flat dimension of any $R$-module is bounded above by
some fixed integer $N$. By \cite[Remark 11.33]{BS}, one can choose $N = 2\dim (R)$.
\end{fact}

\begin{corollary}\label{gab}
Let $(R,\fm, k)$ be a complete local  domain  of prime characteristic. Then $\gd(R^\infty)\leq2\dim (R)+1$.
\end{corollary}

\begin{proof}
 For each positive integer $n$, set $R_n:=\{x\in
R^\infty|x^{p^n}\in R\}$.
This is easy to see that $R_n$ is noetherian and that $R^\infty=\bigcup R_n$. Any countable union of
noetherian ring is
$\aleph_0$-noetherian. We combine Lemma \ref{j} along with Fact \ref{uni}
to observe that $$\pd_{R^\infty}(R^\infty/\fa)\leq\fd_{R^\infty}(R^\infty/\fa)+1\leq2\dim (R)+1.$$
By Auslander's local-global-theorem (please see \cite[Theorem 2.17]{O1}):
$$\gd(R^\infty)=\sup\{\pd(R^\infty/\fa):\fa\lhd R^\infty\}\leq 2\dim( R)+1.$$
\end{proof}

Observation 3.1(xiii) implies that any prime quotient of $R^\infty$ is regular.
Such a thing never happens in the  local algebra. However,
the primeness assumption is important:

\begin{example}
Let $A_0$ be the ring of polynomials with nonnegative
rational exponents in an indeterminant $x$ over a field $\mathbb{F}_{2}$. Let $T$ be
the localization of $A_0$ at $(x^{\alpha}:\alpha>0)$ and set
$A:=T/(x^{\alpha}u : \emph{ u is unit, } \alpha>1)$. In view of
\cite[Page 53]{O1}, $A$ has finite global dimension on maximal
ideals and $\pd(x^{1/2}A)=\infty$.
\end{example}

\begin{theorem}\label{wol}
Let $R$ be a  complete local  domain  of prime characteristic and suppose that its perfect closure is coherent (e.g., $R$ is regular).
The following holds:\begin{enumerate}
\item[i)] If $R$ is not a field, then $\gd(R^{\infty})=\dim (R)+1$.
\item[ii)] Also,  $\Wdim(R^{\infty})=\dim (R)$.
\end{enumerate}
\end{theorem}

\begin{proof}
i)  This is in \cite[Proposition 3.4]{desing}.

ii) One has $\Wdim(R^{\infty})\geq\dim (R)$.  By $i)$, the weak dimension of $R^{\infty}$ is finite.  In view of Theorem \ref{wdim=dim},
  $\Wdim(R^{\infty})\leq\dim (R^{\infty})=\dim (R)$. Thus, $\Wdim(R^{\infty})=\dim (R)$.
\end{proof}

\textrm{Second proof of Theorem \ref{wol}(ii).} Without loss of the generality we assume that $R$ is not a field.
Suppose $\Wdim(R^{\infty})\neq\dim (R)$. Then $\dim (R)+1\leq\Wdim(R^{\infty})\leq\gd(R^{\infty})$. Thus, $R^{\infty}$ is supper regular.
In the light of  Fact \ref{sr} the maximal ideal of $R^{\infty}$ is finitely generated.  But $\fm_{R^{\infty}}^p=\fm_{R^{\infty}}$. By Nakayama's
lemma, $\fm_{R^{\infty}}=0$. This implies that $R^\infty$ is a field. So, $R$ is a field, a contradiction.

\begin{corollary}\label{wdimwe}
Let $R$ be a complete local $\F$-coherent and $S$ be weakly \'{e}tale. Then $\Wdim(S^\infty)\leq\dim R$.
\end{corollary}

In the next  section we reprove (and extend) this by avoiding  Fact \ref{ol}.

\begin{proof}The base change of a flat ring map is flat. This means that $(S\otimes_RR^\infty)\otimes_{R^\infty}(S\otimes_RR^\infty)\to(S\otimes_RR^\infty)$
and $R^\infty\to S\otimes_RR^\infty$ flat.
Thus $R^\infty\to S\otimes_RR^\infty$ is weakly \'{e}tale. By Fact \ref{ol} $\Wdim(S\otimes_RR^\infty)\leq\Wdim(R^\infty)$. Since $\Wdim(R^\infty)=\dim R$, we have $\Wdim(S\otimes_RR^\infty)\leq\dim R$.
 In view of Observation \ref{ob3}(xii), $S^\infty\simeq S\otimes_RR^\infty$. So, $\Wdim(S^\infty)\leq\dim R$.
\end{proof}

\begin{question}
Let $d>0$ be any integer and let $e$ be such that $d+1\leq e\leq 2d+1$.
Is there a $d$-dimensional local ring  $R$  such that
$\gd(R^{\infty})=e$?
\end{question}

\begin{corollary}\label{1}
Let $R$ be a $1$-dimensional complete local domain of prime characteristic. The following are equivalent:
\begin{enumerate}
\item[i)] $R^{\infty}$  is stably coherent.
\item[ii)] $\gd(R^{\infty})=2$.
\item[iii)] $\Wdim (R^{\infty})=1$.
\item[iv)] $R^{\infty}$  is a valuation ring.
\item[v)] $R^{\infty}$  is the perfect closure  of a noetherian regular local ring.
\end{enumerate}
If  $R^{\infty}$ is not coherent, then  $\gd(R^{\infty})=3$.
\end{corollary}

\begin{proof}  Suppose first that $R^{\infty}$ is not coherent.
Any integral domain of global dimension less than $3$ is coherent, see \cite[Theorem 6.3.4]{G}. Thus, $\gd(R^{\infty})>2$.
By Corollary \ref{gab}, $\gd(R^{\infty})\leq2\dim (R)+1=3.$ So, $\gd(R^{\infty})=3$.
\begin{enumerate}
\item[$i)\Rightarrow ii)$] Stably coherent rings are coherent. The claim now follows from Theorem \ref{wol}.
\item[$ii)\Rightarrow iii)$] Note that $0<\Wdim (R^{\infty})\leq\gd(R^{\infty})=2$.
Suppose on the contradiction that  $\Wdim (R^{\infty})\neq1$. Then $\Wdim (R^{\infty})=2$.
In view of \cite[Theorem 6.3.4]{G}, any integral domain of global dimension less than $3$ is coherent.
We conclude from Fact \ref{sr} that $\fm_{R^{\infty}}$  is finitely generated.
But $\fm_{R^{\infty}}$ is not finitely generated, because $\fm_{R^{\infty}}^p=\fm_{R^{\infty}}$. This contradiction shows that
$\Wdim (R^{\infty})=1$.
\item[$iii)\Rightarrow iv)$] In view of \cite[Corollary 4.2.6]{G}
any ring of weak dimension less than $2$ is locally a valuation domain.
Thus, $R^{\infty}$  is a valuation ring.
\item[$iv)\Rightarrow v)$] Any valuation ring is integrally closed. Let $\overline{R}$ be the integral
closure of $R$. Then $$\overline{R}\subset  \overline{R^{\infty}}=R^{\infty}.$$On the other hand
$\overline{R}$ is local, because $R$ is a $1$-dimension complete local ring. Since any $1$-dimensional integrally
closed local domain is discrete valuation ring ($\DVR$), we get that $A:=\overline{R}$ is regular. Clearly,
$A^{\infty}=R^{\infty}$. From this we get the claim.

\item[$v)\Rightarrow i)$] Let $A$ be a  $\DVR$ such that $A^{\infty}=R^{\infty}$. Since $A^{1/p^n}$ is torsion-free
over $A$ it is flat over $A$. Thus, $R^{\infty}$ is a flat filtered limit of noetherian rings. In view of Fact \ref{flatdir}
 $R^{\infty}$  is coherent. Similarly, $R^{\infty}[X]$ is coherent.
\end{enumerate}
\end{proof}

The structure of free resolutions over $R^\infty$ is quite mysterious. However,
 for radical ideals we have the following result:

 \begin{proposition}  \label{str}
Let $(R,\fm)$ be a  Noetherian  local domain  of prime characteristic. The following holds:
\begin{enumerate}
\item[$\mathrm{i)}$] For any radical ideal  $\fa$ of $R^{\infty}$,  the module $R^{\infty}/ \fa$ has a   free resolution of countably generated
free $R^{\infty}$-modules of length  bounded by $2\dim (R)$.

\item[$\mathrm{ii)}$]  For any radical ideal  $\fa$ of $R^{\infty}$, the module $R^{\infty}/ \fa$ has a   flat resolution of countably generated
flat $R^{\infty}$-modules of length  bounded by $\dim (R)$.
\end{enumerate}
\end{proposition}

\begin{proof} $\mathrm{i)}$: By $(x^{\infty})$ we
mean that $(x^{1/p^n}:n\in \mathbb{N}_0)$ where $p:=\Char R$.  Let $d:=\dim (R)$.
  Let $\fa$ be a radical ideal of $R^{\infty}$ and set $\fb:=\fa\cap R$.  Clearly, $\fb$ is radical. By the folklore result of Kronecker,  there is a finite sequence $\underline{\alpha}:=\alpha_1,\ldots,\alpha_{d}$ of elements of $R$ such that $\sqrt{\underline{\alpha}R}=\sqrt{\fb}=\fb$. Suppose $x\in \fa$. Then $x^{p^m}\in R\cap \fa=\fb$ for some integer $m$. It yields that $x^{p^n}=r_1\alpha_1+\cdots+r_d\alpha_d$ for some integer $n$ where $r_i\in R$. By taking $p^n$-th root, $x=r_1^{1/p^n}\alpha_1^{1/p^n}+\cdots+r_d^{1/p^n}\alpha_d^{1/p^n}$. Therefore, $x\in \sum_{i=1}^{d}(\alpha_{i}^{\infty})$, i.e.,
    $\fa\subset\sum_{i=1}^{d}(\alpha_{i}^{\infty})$. The reverse inclusion is trivial. Any radical ideal of $R^{\infty}$ is of the form $\sum_{i=1}^{d}(\alpha_{i}^{\infty})$ for some $\alpha_1,\ldots,\alpha_{d}\in R^{\infty}$. Now we use a trick taken from \cite[Lemma 7.7]{A}.
 Let $F_0$ be a free
$R^\infty$-module with base $\{e_n:n\in \mathbb{N}_0\}$. The assignment
$e_n\mapsto x^{1/p^n}_i$ provides a natural epimorphism
$\varphi:F_0\longrightarrow(x^{\infty}_i)$.
 Let $F_1$ be a free
$R^\infty$-module with base $\{f_n:n\in \mathbb{N}_0\}$. Set $\eta_n:=e_n-x^{\frac{p-1}{p^{n+1}}}_ie_{n+1}$
The assignment
$f_n\mapsto \eta_n$ provides a natural epimorphism
$\varphi:F_1\longrightarrow\ker\varphi$.  Then  a free resolution of $R^{\infty}/(x^{\infty}_i)$  is given by
 the following exact complex:$$\begin{CD}
\textbf{P}_i:
0 @>>> \bigoplus _\mathbb{N}R^{\infty} @>X>> \bigoplus _\mathbb{N}R^{\infty} @>Y>>  R^{\infty} @>>>R^{\infty}/(x^{\infty}_i)@>>>0 \quad(\ast)
\\
\end{CD}
$$
where the matrixes $X$ and $Y$ are defined by:
\begin{equation*}
X := \left(
\begin{array}{cccccc}
1      & 0   & 0 & 0 & 0 &\cdots\\
- x^{\frac{p-1}{p}}_i& 1   & 0 & 0 & 0 &\cdots \\
0 & -x^{\frac{p-1}{p^{2}}} _i  & 1 & 0 & 0 &\cdots \\
0 & 0   & -x^{\frac{p-1}{p^{3}}}_i  & 1 & 0 &\cdots\\
0 & 0   & 0 &  -x^{\frac{p-1}{p^{4}}}_i & 1 & \cdots\\
\vdots & \vdots   &  \vdots& \vdots& \vdots & \ddots\\
\end{array}  \right)
\end{equation*}
and
\begin{equation*}
Y:= \left(
\begin{array}{ccccc}
x _i& x^{1/p} _i& x^{1/p^2}_i & \ldots  \\
\end{array} \right)^t.
\end{equation*}  We apply an easy induction. Set $(\underline{g}^\infty):=\sum_{i=2}^{d}(x_i^{\infty})$.
By the induction hypothesis, $\textbf{P}:=\bigotimes_{i=2}^{d}
\textbf{P}^i$ is  exact. Recall that
$$H_n(\textbf{P}^1\otimes_{R^\infty}\textbf{P})\simeq \Tor^{R^\infty}_n(R^\infty/
(x_1^{\infty}),\frac{R}{(\underline{g}^\infty)}).$$ The ideal
$(x_1^{\infty})$ is  flat. Hence, for each $n>1$ we have
$\Tor^{R^\infty}_n(R^\infty/ (x_1^{\infty}),\frac{R^\infty}{
(\underline{g}^\infty)})=0$. Also,
$$\Tor^{R^\infty}_1(R^{\infty}/ (x_1^{\infty}),\frac{R^\infty}{(\underline{g}^\infty)})
\simeq \frac{(x_1^{\infty})\cap(\underline{g}^\infty)}{
(x_1^{\infty})(\underline{g}^\infty)}=0,$$by Observation \ref{ob3}(v).
This is  clear that $H_0(\textbf{  P}^1\otimes_{R^\infty}\textbf{P})\simeq
R^\infty/(x_1^{\infty})\otimes_{R^\infty}\frac{R^\infty}{(\underline{g}^\infty)}\simeq
R^\infty/ \fa.$
 Thus,
$\bigotimes_{i=1}^d
\textbf{P}^i$ is an explicit
free resolution of $R^\infty/\fa$ of length   $2d$.

$\mathrm{ii)}$:  This is  similar as above.

\end{proof}

\begin{remark}
The above resolution is not necessarily minimal.
\end{remark}

The local assumption is important:

\begin{example}
Let $R$ be  a noetherian regular ring of prime characteristic and infinite global dimension (such a thing exists).
 Then $R^\infty$ is of infinite global dimension.
\end{example}

\begin{proof}Let $n$ be in $\mathbb{N}$. There  is a regular sequence $\underline{x}:=x_1,\ldots,x_n$ over $R$.
This is also is a regular sequence  over $R^\infty$. In particular, $\pd(R^\infty/\underline{x}R^\infty)=n$.
So, $\gd(R^\infty)=\infty$.
\end{proof}

The noetherian assumption is important:

\begin{example}\label{needno}Let  $(A,\fm_A)$ be the localization of $\mathbb{F}_p[X_1,\ldots]$ at $(X_1,\ldots)$. Then $\fd_{A^\infty}(A^\infty/ \fm A^\infty)=\infty$.
\end{example}

\begin{proof}
Let $(A_i,\fm_i):=\mathbb{F}_p[X_1,\ldots,X_i]_{(X_1,\ldots,X_i)}$. In view of Observation \ref{ob3}(ii) $A^\infty$ is flat over $A$. Hence, $$\Tor_i^{A}(A/ \fm_A,A/\fm_A)\otimes_AA^\infty\simeq\Tor_i^{A^\infty}(A^\infty/ \fm A^\infty,A^\infty/\fm A^\infty).$$ To conclude its enough to show that $\Tor_j^{A}(A/ \fm,A/\fm)\neq0$ for all $j$. Let $i$ be a non negative integer. Then, $\Tor_i^{A_i}(A_i/ \fm_i,A_i/\fm_i)\neq0$. By using the rigidity of $\Tor$'s modules (see \cite{au})  we have $\Tor_j^{A_i}(A_i/ \fm_i,A_i/ \fm_i)\neq0$ for any $j\leq i$. Thus, $A_i/ \fm_i\hookrightarrow\Tor_j^{A_i}(A_i/\fm_i,A_i/ \fm_i)$. So, $$0\neq\underset{i}{\varinjlim}A_i/\fm_i\hookrightarrow\underset{i}{\varinjlim}\Tor_j^{A_i}(A_i/\fm_i,A_i/ \fm_i)\simeq\Tor_j^{A}(A/ \fm,A/\fm).$$
\end{proof}

\section{An application: some questions by Shimomoto}

\begin{question}\label{3.1}
(See \cite[Question 2]{ss}) Let $R$ be  local and  $t\in R$  a non-zero divisor. If $R/t R$ is $F$-coherent, then is $R$
also $F$-coherent?
\end{question}

\begin{example} \label{e}  (See \cite[Footnote 24]{BS}) Look at the reduced ring $R:=\mathbb{F}_p[X,Y]/(XY)$. Then $\fd(R^{\infty}/(x))=2\dim (R)$.
In fact they showed that $\Tor^{R^{\infty}}_2(R^{\infty}/(x),R^{\infty}/(y))\neq0$.
\end{example}

Let us give an example of $1$-dimensional ring such that $\gd(R^{\infty})=3$.
To this end, we recall the following beautiful result of Osofsky (see \cite[Proposition 2.36]{O1}):

\begin{fact}\label{osozero}
Let $A$ be a quasilocal ring with  zero-divisors. Then $\gd(A)\geq 3$ and $\Wdim(A)\geq 2$.
\end{fact}

\begin{example}  \label{123}Look at the reduced local ring $R:=\left(\frac{\mathbb{F}_p[X,Y]}{(XY)}\right)_{(x,y)}$. Then
$\gd(R^{\infty})= 3$, $\Wdim(R^{\infty})=2$ and $\dim(R^{\infty})=1$.
\end{example}

\begin{proof}
Indeed, let $A:=\mathbb{F}_p[X,Y]/(XY)$. The lowercase letter will stand for elements of $A$. Set $\fm_{A^\infty}:=(x^{1/p^n},y^{1/p^n}:n\in \mathbb{N})$. It is a maximal ideal of $A$. In view of Example \ref{e},
$$\Tor_2^{A^{\infty}}(A^{\infty}/(x),A^{\infty}/(y))\neq0.$$ Since  $$(x,y)\Tor_2^{A^{\infty}}(A^{\infty}/(x),A^{\infty}/(y))=0,$$
and that $\V((x,y)A^{\infty})\subset\max(A^{\infty})$ we observe that
$$\V\left(\Ann\left(\Tor_2^{A^{\infty}}(A^{\infty}/(x),A^{\infty}/(y))\right)\right)=\{\fm_{A^\infty}\}.$$
Thus, $\Tor_2^{A^{\infty}}(A^{\infty}/(x),A^{\infty}/(y))_{\fm_{A^\infty}}\neq0$. Note that $(A^\infty)_{\fm_{A^\infty}}\simeq R^{\infty}$.
Against $\Ext$, $\Tor$ behaves nicely with respect to the localization over non-noetherian rings, that is $$0\neq\Tor_2^{A^{\infty}}(A^{\infty}/(x),A^{\infty}/(y))_{\fm_{A^\infty}}\simeq \Tor_2^{R^{\infty}}(R^{\infty}/(x),R^{\infty}/(y)).$$
In particular, $\Wdim(R^\infty)\geq 2$ (in fact $\Wdim(R^\infty)=2$). Since $R^\infty$ is quasilocal and has zero-divisor, by Fact \ref{osozero} we observe that
$\gd(R^\infty)\geq 3$. In turns out that $$\gd(R^\infty)= 3=2\dim (R)+1,$$as claimed.
\end{proof}

We also need the following result:

\begin{lemma}(See \cite[Lemma 4.2.3]{G})\label{la}
A quasilocal coherent ring with the property that every principal ideal has finite projective dimension is a domain.
\end{lemma}

\begin{observation}
Let $R$ be local and zero-dimensional. Then $R$ is $F$-coherent.
\end{observation}

\begin{proof}
Set $A:=R_{\red}$. This is a field. By definition, $R^\infty=A^\infty$. In particular, $R^\infty$ is a field and so coherent.
\end{proof}

 The following answers Question \ref{3.1}:

\begin{corollary}
Let $R:=\left(\frac{\mathbb{F}_p[X,Y]}{(XY)}\right)_{(x,y)}$ and set $t:=x-y$. Then $t$ is not zero-divisor, $R/tR$ is  $F$-coherent, and $R$
is not $F$-coherent.
\end{corollary}

\begin{proof}
Clearly,  $t\notin \bigcup_{\fp\in\Ass R} \fp$. We deduce from this that $t$ is not zero-divisor. The following natural isomorphisms $$R/tR\simeq \left(\frac{\mathbb{F}_p[X,Y]/(XY)}{(XY,X-Y)\mathbb{F}_p[X,Y]/(XY)}\right)_{(x,y)}\simeq  \left(\frac{\mathbb{F}_p[X,Y]}{(XY,X-Y)}\right)_{(x,y)}\simeq(\frac{\mathbb{F}_p[Y]}{(Y^2)})_{(y)}$$ implies that $R/tR$ is an $F$-coherent ring (see the above observation).
Combining Lemma \ref{la} and Example \ref{123}, $R$
is not $F$-coherent.
\end{proof}

\begin{remark}\label{newex}In fact, if $R$ is a $1$-dimensional integral domain and $t\neq0$ is any element, then $R/tR$ is $F$-coherent without any regard  with respect
to $F$-coherent property of $R$. In Example \ref{7.3} we will give examples of 1-dimensional  integral domains that are not $F$-coherent.
\end{remark}

Let $I\lhd A$ be  of finite colength.
The \textit{Hilbert-Kunz multiplicity} of an ideal $I$  is defined by $\e_{HK}(I,A):=\lim_{n\to\infty}\frac{\ell \left(\F^n(A/ I)\right)}{p^{n\dim (A)}}.$

\begin{question}(See \cite[Question 3]{ss})
Let $A$ be an $F$-coherent local ring. Is $\e_{HK}(I,A)\in \mathbb{Q}$?
\end{question}

Let $A$ be any commutative ring with an ideal $\frak a$ with a generating set
$\underline{a}:=a_{1}, \ldots, a_{r}$.
 By $\HH_{\underline{a}}^{i}(M)$ we mean the $i$-th cohomology of
\textit{$\check{C}ech$} complex of  a module $M$ with respect  to $\underline{a}$. This is independent of the choose
of the generating set. For simplicity, we denote it by $\HH_{\frak a}^{i}(M)$. This cohomology theory introduced by Grothendieck in SGA 2 \cite{sga2}.

\begin{discussion}
In the case that $\fa\lhd A $ is finitely generated by generating set
$\underline{x}$,  the $\check{C}ech$ grade of
$\fa$ on $M$ is defined by $\inf\{i\in
\mathbb{N}_0 |\HH_{\underline{x}}^{i}(M)\neq0\}.$  Denote it by $\Cgrade_{A}(\fa,M)$. It is easy to observe that $\Cgrade_{A}(\fa,M)=\Kgrade_{A}(\fa,M).$
\end{discussion}
Suppose $(R, \fm)$ is a reduced and of characteristic $p$. Recall that $R$ is $\F$-injective if  $\F: \HH^i_{\fm}(R) \to \HH^i_{\fm}(R)$ is injective.

\begin{lemma}\label{regu}
Let $(R,\fm,k)$ be a noetherian local reduced ring of prime characteristic $p$. If $\underline{x}=x_1,\ldots,x_m$ is a regular sequence over $R$ then
$\underline{x}$ is a regular sequence over $R^\infty$.  Thus, $\depth R\leq\Kdepth_{R^\infty} (R^\infty)$. The equality holds if $R$ is $F$-injective or Cohen-Macaulay.
\end{lemma}

\begin{proof}
Let $r\in R^\infty$ be such that $rx_i\in(x_1,\ldots,x_{i-1})R^\infty$. There are $r_j\in R^\infty$ such that
 $rx_i=\sum_{j=1}^{i-1}r_jx_j.$ Recall that $R^\infty=\bigcup R^{1/p^n}$. In particular, there is an $n$ such that
$r$ and all of $r_j$ are in $R^{1/p^n}$. Taking $p^n$-th power we have $$r^{p^n}x_i^{p^n}=\sum_{j=1}^{i-1}r_j^{p^n}x_j^{p^n}\quad(\ast)$$Recall from \cite[Theorem 16.1]{Mat} that $\underline{x}^{p^n}:=x_1^{p^n},\ldots,x_m^{p^n}$ is a regular sequence over $R$. Since every thing in $(\ast)$ is in $R$,
we have $r^{p^n}\in(x_1^{p^n},\ldots,x_{i-1}^{p^n})R$. Taking $p^n$-th root, we have $$r\in(x_1,\ldots,x_{i-1})R^{1/p^n}\subset (x_1,\ldots,x_{i-1})R^{\infty}.$$ Thus,
$\underline{x}$ is a regular sequence over $R^\infty$.  Therefore, $\depth_R(R)\leq\depth_{R^\infty}(R^\infty)\leq\Kdepth _{R^\infty}(R^\infty).$

Suppose now that $R$ is $F$-injective. Since $R$ is $F$-injective, then $\HH^i_{\fm}(R)\stackrel{\F}\lo\HH^i_{\fm}(R)$ is injective by definition.
Since, the local
    cohomology commutes with direct limit (via Frobenius), $\HH^i_{\fm}(R)\lo\HH^i_{\fm}(R^\infty)$ is injective. Hence, $\Cgrade_{R}(\fm,R)\geq \Cgrade_{R}(\fm,R^\infty)$  $(\dagger)$. Note that Koszul grade is the same as of $\check{C}ech$ grade.
 Thus, \[\begin{array}{ll}
\depth_R(R)&\geq \Kdepth_R(R^\infty)\ \  \ \ \ \ \ \ \ \ \ \ \ \ \ \ \ \ \  \ \ \ \ \ \ \ \ \  \ \ \ \ \ \ \  \ by \ \ (\dagger)   \\
&\stackrel{}=\Kgrade_R(\fm,R^\infty)\ \ \ \ \ \ \ \ \ \ \ \ \ \ \ \ \ \  \ \ \ \ \ \ \ \ \ \ \textit{ by Definition  }  \\
&\stackrel{}=\Kgrade_{R^\infty}(\fm R^\infty,R^\infty) \ \ \ \ \ \ \ \ \ \ \ \ \ \ \ \ \ \ \ \ \  \textit{by Fact }\ref{cp}(iii)\\
&\stackrel{}=\Kdepth_{R^\infty}(R^\infty) \ \ \ \ \ \  \ \ \ \ \ \ \ \ \ \ \ \ \ \ \ \ \ \ \ \ \ \ \ \ \ \ \ \textit{by Fact }\ref{cp}(i)\\&\geq\depth_R(R)\ \ \ \ \ \  \ \ \ \ \ \ \ \ \ \ \ \ \ \ \ \ \ \ \ \ \ \ \ \ \ \ \ \ \ \ \  \ \ \ \ \textit{by the first part }
\end{array}\]Thus, $\depth R=\Kdepth (R^\infty)$, as claimed.

Finally, suppose that $R$ is Cohen-Macaulay. In view of Fact A) in Theorem \ref{wdim=dim}, $\Kdepth(R^\infty)\leq\dim (R^\infty).$
By the first part, we have $$\dim (R^\infty)=\dim R=\depth R\leq \Kdepth(R^\infty)\leq\dim (R^\infty).$$
The proof is now complete.
\end{proof}

\begin{remark}i) In the above lemma, the $\F$-injectivity assumption is important. There is a non-Cohen-Macaulay
$\F$-coherent ring $R$, see \cite[Example 3.6]{ss}.
By \cite[Theorem 3.11]{ss} $R^\infty$  is big Cohen-Macaulay. Thus, $\depth R<\dim R=\depth (R^\infty)=\Kdepth (R^\infty)$.

ii) Let $R$ be a reduced $\F$-coherent local ring which
is a residue class ring of a Gorenstein local ring. It may worth to recall from \cite[Corollary 3.16]{ss} that $\F$-injectivity
implies Cohen-Macaulayness and $\F$-rationality.
\end{remark}
We need the following result of Seibert:

\begin{fact} (See \cite[Proposition 2(b)]{sib})\label{sib}
Let $\mathcal{C}_{M}$ be the class of all finite $R$-modules $P$ such that the length of
$M\otimes P$ is finite. If $\pd M<\infty$ and $N\in \mathcal{C}_{M}$ then there are certain  $b_i\in \mathbb{Q}$ such that
$$\sum_{i=0}^{\dim R}(-1)^i\ell(\Tor^i(\F^n(M),N))=\sum_{i=0}^{\dim N} b_ip^{in}.$$
\end{fact}
Also, $R$ is called $\F$-finite if $R$ viewed as an $R$-module
via $\F$ is finite. For example, every ring which is a localization of an affine
algebra over a perfect field and every complete local ring with perfect
residue  field is $\F$-finite.

\begin{proposition}
Let $(R,\fm,k)$ be an $\F$-finite and  $\F$-coherent domain  and $J\lhd R$ be  of finite colength. If $R$ is  Cohen-Macaulay, then $\e_{HK}(J,R)$ is rational.
\end{proposition}

\begin{proof}
We know that  $\pd(R^\infty/J R^\infty)<\infty$. Note that $R^\infty/J R^\infty$ is finitely presented. Since $R^\infty$ is coherent, $R^\infty/J R^\infty$ has finite free resolution by finitely generated free modules:$$\xymatrix{0\ar[r]&F_{N}\ar[r]&\ldots\ar[r]&F_{j+1}\ar[r]^{f_j}&F_{j}\ar[r]&\ldots\ar[r]&F_{0}\ar[r]&R^\infty/J R^\infty\ar[r]&0,}$$
There is an index $i\in \mathbb{N}$ such that all of components of  $\{f_j\}$ are in $(R_i,\fm_{R_i})$. Let $F_{j}(i)$ be the
free $R_i$-module with the same rank as $F_{j}$. Consider $f_j$ as a matrix over $R_i$, and denote it by $f_j(i)$.
Recall that $\fm$ is finitely generated. Look at the following
complex of finite free  modules:$$\xymatrix{(\ast)\quad0\ar[r]&F_{N}(i)\ar[r]&\ldots\ar[r]&F_{j+1}(i)\ar[r]^{f_j(i)}&F_{j}(i)
\ar[r]&\ldots\ar[r]&F_{0}(i)\ar[r]&R_i/J R_i\ar[r]&0.} $$Note that  $(\ast)$ is a   complex.
We are going to show that $(\ast)$ is exact.  Let $r_j$ be the expected rank of $f_j$.  Recall that $I_{t}(f_j(i))$ is the ideal
generated by $t\times t$ minors of $f_j(i)$. Clearly,  $r_j$ is
the expected rank of $f_j(i)$.

The localization of Cohen-Macaulay is
again Cohen-Macaulay. For each $\fp \in \Spec(R^\infty)$, we set $P:=\fp\cap R_i$.   Then
\[\begin{array}{ll}
j&\leq\Kgrade_{(R^\infty)}(I_{r_j}(f_j), R^\infty)\ \ \ \ \ \ \ \ \  \ \ \ \ \ \ \ \ \ \ \ \ \ \ \ \ \ \ \ \ \ \ \  \ \ \ \ \ \ \ \ \ \ \  \ \ \ \ \ \ \ \ \  \ \ \ \ \ \  \ \ \ \ \ \ \ \ \  \ \
Fact\ \ \ref{cp}(v)\\
&=\inf\{\Kdepth_{(R^\infty)_{\fp}}((R^\infty)_{\fp}):\fp\in V(I_{r_j}(f_j))\}\ \ \ \ \ \ \ \ \ \ \ \ \ \ \ \ \ \ \ \ \ \ \  \ \  \ \ \ \ \ \  \ \ \ \ \ \ \ \  Fact \ \ \ref{cp}(iv)\\
&=\inf\{\Kdepth_{((R_i)_P)^\infty}((R_i)_P)^\infty:P\in V(I_{r_j}(f_j)(i))\} \ \ \ \ \ \  \ \ \ \ \ \ \ \ \ \  \ \ \ \ \ \  \ \ \ \  \ \ \textit{ Observation \ref{ob3}(xi)}\\
&=\inf\{\Kdepth_{(R_i)_P}((R_i)_P):P\in V(I_{r_j}(f_j(i))\}\ \   \ \ \ \ \ \ \ \ \ \ \ \ \ \ \ \ \ \ \ \ \ \ \ \ \ \ \ \ \ \ \ \ \ \ \ \    \textit{ Lemma  }   \ref{regu}\\
&\stackrel{}=\Kgrade_{R_i}(I_{r_j}(f_j(i),R_i) \ \ \ \ \ \ \ \ \ \ \ \ \   \ \ \ \ \ \ \   \ \ \ \ \ \ \   \ \ \ \ \  \ \ \ \ \ \ \ \ \ \ \ \ \ \ \ \ \ \ \ \ \ \ \ \ \ \ \ \ \ \ \ \ \ \ \ \ \ \ \ \ Fact \ \ \ref{cp}(iv)  \\
\end{array}\] Again, by  applying Fact \ref{cp}(v), $(\ast)$ is exact.  Thus, $R_i/J R_i$ is of finite projective dimension.
By Fact \ref{sib} $\e_{HK}(J R_i,R_i)\in \mathbb{Q}$. Since  $R$ is $F$-finite,  $R_i$ is finitely generated as a module over $R$.
Combine this along with  \cite[Theorem 2.7]{wa} we get that
$$\e_{HK}(J,R)=\frac{\e_{HK}(J R_i,R_i)[R_i/ \fm_{R_i}:k]}{[Q(R_i):Q(R)]}\in \mathbb{Q},$$where $Q(-)$ stands for the fraction field of an integral domain. The proof is now complete.
\end{proof}

The following extends \cite[Proposition 3.19]{ss} by Shimomoto where he worked with $R\to R^{\hens}$. This may answer
\cite[Question 1]{ss} where he asked conditions on the fibers of a flat extension to ascend the $\F$-coherent property:

\begin{proposition}\label{cohwe}
Let $A$ be a  local $\F$-coherent and $B$ be weakly \'{e}tale over $A$. Then $B$ is $\F$-coherent.
\end{proposition}

\begin{proof}
Recall that $A^\infty\to B\otimes_AA^\infty$ is weakly \'{e}tale. By definition $A^\infty$ is coherent as
an $A^\infty$-module. By \cite[Proposition]{Ol}, $B\otimes_AA^\infty$ is coherent as
an $B\otimes_AA^\infty$-module. In view of Observation \ref{ob3}(xii), $B^\infty\simeq B\otimes_AA^\infty$.
Thus, $B^\infty$ coherent. I.e., $B$ is $\F$-coherent.
\end{proof}

From Proposition \ref{cohwe} one can recover Corollary \ref{wdimwe}.

\section{homological dimension over  $R^\textbf{+}$}

We start by recalling some historical remarks:

\begin{discussion}\label{hist}
Hochster
proved  that $\fd_{R^\textbf{+}}(R^\textbf{+}/ \fm_{R^\textbf{+}})\leq\dim (R)$, when $R$ is  henselian
and has residue prime characteristic, see \cite[Proposition 2.15]{H1}.
Hochster and Aberbach extended this  by showing  that the
flat dimension of any radical of a finitely generated ideal has finite flat dimension, please see \cite[Theorem 3.1]{ab}.
Also, recall from \cite[Theorem 1.1(i)]{A} that
$\pd_{R^\textbf{+}}(R^\textbf{+}/ \fm_{R^\textbf{+}})\leq2\dim (R)$.
 \end{discussion} This lead us to ask:

\begin{conjecture}
 Let $(R,\fm)$ be a  complete local domain of  prime characteristic. Then $R^\textbf{+}$ is regular.
\end{conjecture}

In this section we assume the generalized continuum hypothesis that is $2^{\aleph_n}=\aleph_{n+1}$.
We will use this only for computing projective dimension (but not for flat dimension). Is it really needed?

\begin{lemma}\label{j1} Let $(R,\fm, k)$ be a complete local  domain  of prime characteristic. If $k$ is countable, then
$R^\textbf{+}$ is of $\aleph_1$ cardinality.
\end{lemma}

\begin{proof}
By Cohen's Structure theorem, $R$ is  a module-finite extension of a complete regular local ring $A$. It is not difficult to see that
$A^\textbf{+}=R^\textbf{+}$ and that the residue field of $A$ is countable. Then without loss of the generality we may assume that
$R$ is complete. Again by Cohen's Structure theorem $R$ is of the form $R:=k[[x_1,\ldots,x_d]]$. Any element of $R$
is a formal power series with coefficient taken from $k$.  Thus the cardinality of $R$ is  the cardinality of $\prod_\mathbb{N} k$.
By $|-|$ we mean the cardinality of a set. Hence, $|R|=\aleph_1$. Denote the fraction field of $R$ by $Q$. Since $Q:=\{r/s: r\in R, s\neq 0\}$, we observe that $Q$ is the same cardinality as of $R$. So, $|Q|=\aleph_1$. Let $\overline{Q}$ be the algebraic closure of $Q$. We are going to show that $|\overline{Q}|=\aleph_1$.
By definition, $\overline{Q}$ determines by the root of polynomial with coefficient in $Q$. Note that $|Q[X]|=\aleph_1$. From this we deduce that
$|\overline{Q}|=\aleph_1$. Since $R^\textbf{+}\subset\overline{Q}$, we get the claim.
\end{proof}

\begin{proposition}\label{cc}
Let $(R,\fm, k)$ be a complete local  domain  of prime characteristic. The following holds:
\begin{enumerate}
\item[i)] Any $R^\textbf{+}$-module has a finite flat dimension. In fact, $\Wdim(R^\textbf{+})\leq2\dim (R)$.
\item[ii)] If $k$ is countable, then $R^\textbf{+}$ is regular. In fact, $\gd(R^\textbf{+})\leq2\dim (R)+1$.
 \end{enumerate}
\end{proposition}

\begin{proof} i) Let $\fa$ be an ideal of $R^\textbf{+}$.
Recall that $R^\textbf{+}=\bigcup R_\gamma$, where $R_\gamma$ is a module-finite extension of $R$.
Without loss of the generality we may assume that
 $R_\gamma$ is complete and local. Indeed,
 let $\overline{R_\gamma}$ be the integral closure of $R_\gamma$ in its
field of fractions. Recall that the
integral closure of a complete local domain in its field of
fractions is  local and complete. Thus $\overline{R_\gamma}$ is a
 complete local  normal domain. To conclude, it  remains
to recall that $\overline{R_{\gamma}}\subseteq R^\textbf{+}$. Write
 $R^\textbf{+}=\bigcup R_\gamma$, where $R_\gamma$ is a
complete  local domain of dimension equal to $d:=\dim (R)$. Note that
$R^\textbf{+}\simeq{\varinjlim}R_\gamma^\infty.$
Let $\fa$ be an ideal of $R^\textbf{+}$.
Set $\fa_\gamma:=\fa \cap R_\gamma^\infty$.  One has  ${\varinjlim}R_\gamma^\infty/ \fa_\gamma  \simeq R^\textbf{+}/\fa.$ Recall from  \cite[VI, Exercise 17]{CE} that $$\underset{i}{\varinjlim}\Tor_j^{R_\gamma^\infty}(R_\gamma^\infty/ \fa_\gamma,-)\simeq\Tor_j^{R^\textbf{+}}(R^\textbf{+}/ \fa,-).$$
Let $j>2\dim (R)$. By applying Fact \ref{uni}, we conclude that $\Tor_j^{R^\textbf{+}}(R^\textbf{+}/ \fa,-)=0.$
Since $R$ is local, its dimension is finite. Thus, $$\fd(R^\textbf{+}/ \fa)<2\dim (R)+1<\infty.$$
Due to $(2.1.1)$ any module has finite flat dimension.

ii)  Let $\fa$ be an ideal of $R^\textbf{+}$. In view of Lemma \ref{j1}, $ R^\textbf{+}$ is $\aleph_0$-noetherian. Keep part i) in mind. In view of Lemma \ref{j} we see that $\pd_{R^\textbf{+}}(R^\textbf{+}/\fa)\leq\fd_{R^\textbf{+}}(R^\textbf{+}/\fa)+1$.
Thus, $$\pd_{R^\textbf{+}}(R^\textbf{+}/\fa)\leq\fd_{R^\textbf{+}}(R^\textbf{+}/\fa)+2\leq2\dim (R)+1<\infty.$$By Auslander's local-global-theorem (please see \cite[Theorem 2.17]{O1}):
$$\gd(R^\textbf{+})=\sup\{\pd(R^\textbf{+}/\fa):\fa\lhd R^\infty\}\leq 2\dim( R)+1.$$
\end{proof}

The above bound may not be sharp:

\begin{example}
If  $R$ is one-dimensional local and complete,  then $$\gd
(R^\textbf{+})=\pd_{R^\textbf{+}}(R^\textbf{+}/ \fm_{R^\textbf{+}})=2=2\dim( R)$$ and $\Wdim( R^\textbf{+})=\dim( R)=1$.
\end{example}

\begin{proof}
The first claim is in \cite[Theorem1.1(ii)]{A} where under this assumption we observed that $R^\textbf{+}$ is a valuation domain. So,
$\Wdim( R^\textbf{+})=\dim( R)=1$.
\end{proof}

Since  \cite[Lemma 4.1]{A}  is  true via the henselian assumption, we take this opportunity to state its corrected
version:

\begin{lemma}\label{art}(Artin)
Let $(R,\fm, k)$ be a henselian (e.g. complete) local  domain  of prime characteristic. Then $R^\textbf{+}$ is quasilocal.
\end{lemma}

The henselian assumption is important. Let us  analyze this by the help of an explicit example.

\begin{example}(Epstein)
It may be $R$ is local but $R^\textbf{+}$ is not quasilocal.

i) Let $R:=\mathbb{Z}_{(p)}$, where $p$ is an odd prime number.  Let $S:=R[\sqrt{p+1}]$.  This ring is integral over $R$. In particular, $\dim( S)=\dim( R)=1$  and that $R^\textbf{+} = S^\textbf{+}$. Look at $\alpha := \sqrt{p+1}-1$ and  $ \beta := \alpha+2$.  Of course $\beta - \alpha = 2$ which is a unit of $R$. Hence $\beta - \alpha$ is a unit of $S$.  This implies that $\beta,$ and $ \alpha$ cannot live together in any prime ideal.  Neither $\alpha$ nor $\beta $ is a unit of $S$.
Indeed, note that $\alpha \beta = p$. Thus, the inverse of $\alpha$ in the real numbers is $\beta / p =  \frac{\sqrt{p+1}}{p} + (1/p)$, which is clearly not in $S$.  Similarly, the inverse of $\beta$ in the real numbers is $\alpha/p$, also not in $S$.
Therefore, there must be some prime ideal $\fp\lhd S$ that contains  $\alpha$, and a different prime ideal $\fq$ of $S$ that contains  $\beta$.  Now, note that $R^\textbf{+} = S^\textbf{+}$.  By the lying-over property of integrality, there must be prime ideals $\fp'$ and $\fq'$ of $R^\textbf{+}$ that contract to $\fp$, and $\fq$ respectively. Both  of $\fp'$ and $\fq'$ must be maximal ideals of $R^\textbf{+}$, because $\dim (R^\textbf{+})=1$.  So $R^\textbf{+} $ is not quasilocal.

ii) One may ask an example in the prime characteristic case. Let $R :=  \mathbb{F}_p[x]_{(x)}$. After replacing $p$ in the first item i) by $x$,
it is routine to see that  $R^\textbf{+}$ is not quasilocal.
\end{example}

Both of the above examples are normal.  If the base ring is normal, the following holds:

\begin{fact}(See \cite[Proposition 1.4]{Ar})
A normal integral domain $A$ is henselian if and
only if $A^\textbf{+}$  is quasilocal.
\end{fact}

\begin{corollary}
Let $(R,\fm, k)$ be a complete local  domain  of prime characteristic. If $k$ is of cardinality $\aleph_n$ (e.g. $k$ is countable), then any ideal of $R^\textbf{+}$ has finite free resolution (by
  not necessarily finitely generated  free modules).
\end{corollary}

\begin{proof}
In view of Lemma \ref{art},  $R^\textbf{+}$ is quasilocal. Since $k$ is of cardinality $\aleph_n$, $R^\textbf{+}$ is  $\aleph_{n+1}$-noetherian, see the proof of Lemma \ref{j1}. Let $\fa$ be an ideal of  $R^\textbf{+}$. By the proof of Proposition \ref{cc},$$\pd_{R^\textbf{+}}(R^\textbf{+}/\fa)\leq\fd_{R^\textbf{+}}(R^\textbf{+}/\fa)+(n+2)\leq2\dim (R)+(n+2)<\infty.$$By a famous result of Kaplansky,
  any projective module  over a quasilocal  ring is free. The proof is now complete.
\end{proof}

\section{Failure of coherence $R^\textbf{+}$ in mixed characteristic }

Our aim is to understand the higher-dimensional version of the following observation.

\begin{observation}Let $(R,\fm,k)$  be a complete regular local ring.
If $\dim (R)= 1$, then $R^\textbf{+}$ is  a filtered colimit of finitely presented flat $R$-algebras. The transition maps are  flat.
In particular,  $R^\textbf{+}$ is coherent.
\end{observation}

\begin{proof}
Recall that $R^\textbf{+}=\bigcup R_\gamma$, where $R_\gamma$ is a module-finite extension of $R$.
Without loss of the generality we may assume that
 $R_\gamma$ is normal and local. In particular,  $R_\gamma$  is $\DVR$. Since torsion-free
modules over  a $\DVR$ are flat, we get the first claim. The particular case follows by Fact \ref{flatdir}.
\end{proof}

\begin{lemma}\label{kunz}(See \cite[Theorem 23.1]{Mat})
Let $\varphi$ be a local map from a  regular local ring $(A,\fm)$ to a  Cohen-Macaulay local ring $(B,\fn)$.
Suppose $\dim (A)+\dim (B/\fm B)=\dim (B)$. Then $\varphi$ is flat.
\end{lemma}

\begin{proposition}Let $(R,\fm,k)$  be a complete regular local ring.
If $\dim (R)= 2$, then $R^\textbf{+}$ is  a filtered colimit of finitely presented flat $R$-subalgebras. The transition maps are not flat
provided $k$ is of positive transcendence degree over $\mathbb{F}_p$.
\end{proposition}

\begin{proof}
Recall that $R^\textbf{+}=\bigcup R_\gamma$, where $R_\gamma$ is a module-finite extension of $R$.
Without loss of the generality we may assume that
 $R_\gamma$ is complete and local. Indeed,
 let $\overline{R_\gamma}$ be the integral closure of $R_\gamma$ in its
field of fractions. Recall that the
integral closure of a complete local domain in its field of
fractions is  local and complete. Thus $\overline{R_\gamma}$ is a
 complete local  normal domain. For simplicity,  $R_\gamma:=\overline{R_\gamma}$. We are going to show that $(R,\fm)\to (R_\gamma,\fm_\gamma)$ is flat.
Since $R_\gamma$ is normal and $2$-dimensional,   Serre's characterization of normality
implies that $R_\gamma$ is   Cohen-Macaulay.
Note that $R$ is regular,  $R_\gamma$ is   Cohen-Macaulay  and that $$\dim (R_\gamma)=\dim (R)+ \dim(\frac{ R_\gamma}{\fm R_\gamma}).$$
 In view of Lemma \ref{kunz}, $R\to R_\gamma$ is flat. This finishes the proof of first claim.

Suppose $t:=\trdeg(k/\mathbb{F}_p)>0$ and suppose on the contrary that $R_\gamma\to R_\beta$ is flat for all $\gamma<\beta$. Then $R^\textbf{+}$ is  a flat  filtered  limit of noetherian rings. In the light of Fact \ref{flatdir}, this implies that $R^\textbf{+}$ is coherent. It is shown in  \cite[Theorem 4.9]{ab} that  $R^\textbf{+}$ is not coherent (here, we need $t>0$). This contradiction shows that $R_\gamma\to R_\beta$ is not flat for all cofinal pair $\gamma<\beta$.
\end{proof}

\begin{corollary}Let $(R,\fm,k)$  be a complete regular local ring.
If $\dim (R)= 2$ and of prime characteristic, then $R^\textbf{+}$ is  a filtered colimit of finitely presented Cohen-Macaulay $R$-subalgebras. Also,
 $R^\textbf{+}$ is not a filtered colimit of its regular local subalgebras over $R$ provided $\trdeg(k/\mathbb{F}_p)>0$.
\end{corollary}

\begin{corollary}Let $(R,\fm,k)$  be a complete regular local ring. Then
 $R^\textbf{+}$ is   a filtered colimit of finitely presented Cohen-Macaulay $R$-subalgebras if and only if  $R^\textbf{+}$ is filtered colimit of finitely presented flat $R$-subalgebras.
\end{corollary}

\begin{proof}
This follows by Lemma \ref{kunz}.
\end{proof}

\begin{theorem}
Let $R$ be a complete domain having mixed characteristic $p$. If
$\dim (R )> 3$ then $R^\textbf{+}$ is not coherent.
\end{theorem}

What can say when $\dim (R)=3$?\footnote{If  the answer is positive, then one can show that $R^\textbf{+}$ is balanced big Cohen-Macaulay. The last one  is an open question.}

\begin{proof}
Let $d:=\dim (R)$ and $p:=\Char (R/ \fm)>0$.
In the light of  \cite{an} the direct summand conjecture is now a beautiful theorem. By Fact \ref{ab} we know that
$\fd(R^\textbf{+}/\fm_{R^\textbf{+}})=d.$ Suppose on the contradiction that $R^\textbf{+}$ is coherent. Let $M$ be a finitely generated $R^\textbf{+}$-module.
In view of \cite[Corollary 2.5.10]{G} $\pd_{R^\textbf{+}}(M)\leq n$ if and only if
$\Tor_{n+1}^{R^\textbf{+}}(M,R^\textbf{+}/\fm_{R^\textbf{+}})=0$.
From this we observe that $\Wdim(R^\textbf{+})=d<\infty$. In view of Theorem \ref{wdim=dim} $$\dim (R)=\Wdim (R^\textbf{+})=\Kdepth( R^\textbf{+}).$$ Let  $\underline{x}:=p,x_2,\ldots,x_d$ be a system of parameters for $R$. Note that
$\rad(\underline{x})=\fm_{R^\textbf{+}}.$ By
 Corollary \ref{should}, $\underline{x}$ is a regular sequence over $R^\textbf{+}$. Recall  from \cite[Proposition 3.6]{ab} that  $R^\textbf{+}$ is not a balanced big Cohen-Macaulay algebra for $R$ (here, we need $d>3$).
This contradiction shows that  $R^\textbf{+}$ is not coherent.
\end{proof}

Also, the following left unsolvable:

\begin{question}
Suppose $k$ is the algebraic closure of $\mathbb{F}_p$. Let $(R, \fm, k)$ be a 2-dimensional complete
normal domain. Is $R^\textbf{+}$ coherent? Specially: Is $I^{\ast}=I^\textbf{+}$? (The first one is tight closure and the second one is plus closure. Also, the first question implies the second.)
\end{question}
Brenner proved that $I^{\ast}=I^\textbf{+}$ where $R$ is a 2-dimensional standard graded ring over algebraic closure of $\mathbb{F}_p$ and for homogeneous ideal $I$. This uses
the theory of vector bundles, see \cite{Br}. His assumption over $R_0$ is really important.

\section{$\gd(R^{\infty})$ depends on the characteristic}

Let $R$ be one-dimensional local and complete. Then $\gd
R^\textbf{+}=2$ and $\Wdim R^\textbf{+}=1$ without any regards with respect to the characteristic.
Here we show $\gd(R^{\infty})$ and $\Wdim(R^{\infty})$ depend to the characteristic.

\begin{lemma}Let $p\neq 2$ be a prime integer.
Let $R:=\mathbb{F}_p[x,y]/(y^2-x^3-x^2)$. The following holds:
\begin{enumerate}
\item[i)] $\gd(R^{\infty})= 3$, and
\item[ii)] $\Wdim(R^{\infty})=2$.
 \end{enumerate}
\end{lemma}

\begin{proof} Fist we claim that $R^\infty$ is not coherent. We follow some lines from \cite{AS}.
Note that
$R\simeq\mathbb{F}_p[t,t\sqrt{t+1}].$ The normalization of $R$ in its field of fractions is $\mathbb{F}_p[t,\sqrt{t+1}]$.  Suppose on the contrary that $R^\infty$ is coherent. Over $1$-dimensional reduced rings, this is equivalent with the property  that $\overline{R}/R$ is purely inseparable, see \cite[Corollary 3.9]{ss}. Since
$(\sqrt{t+1})^{p^n}\notin R$ for each $n$, the extension $\overline{R}/R$
is not  purely inseparable. This contradiction shows that $R^{\infty}$ is not coherent.

i) Recall from
\cite[Theorem 6.3.4]{G}
that a domain of global dimension less than $3$ is coherent. By the first paragraph, $R^{\infty}$ is not coherent. So, $\gd(R^{\infty})\geq 3$.
Due to \cite[Remark 11.33]{BS} $\gd(R^{\infty})\leq2\dim (R)+1= 3$. Therefore, $\gd(R^{\infty})= 3$.

ii) By part $i)$ we have $\gd(R^{\infty})=3$. Recall from   Lemma \ref{j} that $$3=\gd(R^{\infty})\leq 1+\Wdim(R^{\infty}).$$  So, $\Wdim(R^{\infty})\geq 2$.
Due to \cite[Remark 11.33]{BS} $\Wdim(R^{\infty})\leq 2\dim R=2$. Thus   $2\leq\Wdim(R^{\infty})\leq 2 $  as claimed.
\end{proof}

\begin{lemma}
Look at the integral domain $R:=\mathbb{F}_2[t,t\sqrt{t+1}]$. Then
$\gd(R^{\infty})=2$  and $\Wdim(R^{\infty})=1$.
\end{lemma}

\begin{proof}
Look at the regular ring $A:=\mathbb{F}_2[t]$. It is easy to see that $R^{\infty}=A^{\infty}$.
Since $A$ is regular and of dimension one, it is easy to observe that $\gd(R^{\infty})=2$ and $\Wdim(R^{\infty})=1$.
\end{proof}

\begin{example}\label{7.3}
Let $R:=\mathbb{F}_p[t,t\sqrt{t+1}]$. Then
\begin{equation*}
\gd(R^{\infty})= \left\{
\begin{array}{rl}
3 & \  \   \   \   \   \ \  \   \   \   \   \ \text{if } p\neq 2\\
2 & \  \   \   \   \   \ \  \   \   \   \   \ \text{if } p=2
\end{array} \right.
\end{equation*}
Also,
\begin{equation*}
\Wdim(R^{\infty})= \left\{
\begin{array}{rl}
2 & \  \   \   \   \   \ \  \   \   \   \   \ \text{if } p\neq 2\\
1 & \  \   \   \   \   \ \  \   \   \   \   \ \text{if } p=2
\end{array} \right.
\end{equation*}
\end{example}

\begin{proof}
This is the combination of the above lemmas.
\end{proof}

\section{Revisiting the miraculous formula}

We give a  simple proof of the following funny fact (see \cite[Lemma 3.16]{BS}).

\begin{fact}(Bhatt and Scholze) Let
$B\leftarrow A\rightarrow C $ be a diagram of perfect rings of prime characteristic. Then $\Tor^{A}_i (B,
C)=0$ for all $i\geq1$.
\end{fact}

\begin{proof}
Let $\Gamma$  be a generating set of $B$ as an $A$-algebra.  Look at perfection of the polynomial ring $D:=(A[x:x\in \Gamma])^\infty$. In particular,
the map $A\to B$ may view  as   $A\to D\twoheadrightarrow B$.
There ia an spectral sequence $$\Tor^{D}_{i}(\Tor^{A}_j(C,D),B)\Rightarrow \Tor^{A}_{i+j}(B,C).$$
Since $D$ is free as an $A$-algebra, the spectral sequence collapses.  Set $E:=C\otimes_AD$  which is a perfect algebra by Observation \ref{ob3}(x). Thus, $\Tor^{D}_{i}(E,B)\simeq\Tor^{A}_{i}(B,C)$. We apply the replacement $E\mapsto C$ and $D\mapsto A$. Hence,
we may assume that $A \twoheadrightarrow B$ is surjective.  Perfect rings are reduced. So, $B\simeq A/\fa$ for some radical ideal $\fa$. Write
 $A=\bigcup R_\gamma$, where $R_\gamma$ is a noetherian subring of $A$. Taking perfection, we have
$A\simeq{\varinjlim}R_\gamma^\infty.$ Since $\Tor$ modules behave nicely with respect to direct limits we may assume that $A=R^\infty$ for some  noetherian ring $R$. Any radical ideal $\fa$ of $R^\infty$
is of the form $\sum_{i=1}^m (x_i^\infty)$ for some $m\in \mathbb{N}$ where $(x^\infty):=(x^{1/p^n}:n\in \mathbb{N})$. Similarly, $C=A/\fb$
for some radical ideal $\fb$.

First, we deal with the case that $m=1$.  Clearly, $\fd(\frac{R^\infty}{(x^\infty)})\leq1$. Hence, the only crucial $\Tor_i$ is $\Tor_1$:$$\Tor^{R^\infty}_1 (\frac{R^\infty}{(x^\infty)}, R^\infty/ \fb)\simeq\frac{(x^\infty)\cap \frak b}{(x^\infty)\fb}$$which is zero by Observation \ref{ob3}(v).
Set $\fc:=\sum_{i=2}^m (x_i^{1/p^n}:n\in \mathbb{N})$. By induction we have
$$\Tor^{R^\infty}_i (\frac{R^\infty}{\fc},R^\infty/ \fb)=0\quad\forall i>0\quad(\ast)$$There is a change of rings $R^\infty\twoheadrightarrow\frac{R^\infty}{\fc}=:S$. Note that $S$ is the perfect closure of $R/ (\fc\cap R)$. Look at the $S$-module $ R^\infty/ \fa$ and  $R^\infty$-module $ R^\infty/ \fb$, there ia an spectral sequence $$\Tor^{\frac{R^\infty}{\fc}}_{i}(\Tor^{R^\infty}_j(\frac{R^\infty}{\fb},{\frac{R^\infty}{\fc}}),R^\infty/ \fa)\Rightarrow \Tor^{R^\infty}_{i+j}(\frac{R^\infty}{\fb},R^\infty/ \fa).$$
By $(\ast)$ the spectral sequence collapses. Combine this along with the case $m=1$ and the natural isomorphism $S/ (x_1^\infty)\simeq R^\infty/ \fa$ we have$$\Tor^{R^\infty}_{i+j}(\frac{R^\infty}{\fa},R^\infty/ \fb)\simeq\Tor^{\frac{R^\infty}{\fc}}_{i+j}(\Tor^{R^\infty}_0(\frac{R^\infty}{\fb},{\frac{R^\infty}{\fc}}),R^\infty/ \fa)\simeq\Tor^{S}_{i+j}(\frac{S}{\fb},S/ (x_1^\infty))=0.$$The proof is now complete.
\end{proof}

Wodzicki has constructed an example of a ring $A$ such that  $\Jac(A)\neq0$ and  $\Tor^{A}_{\geq1} (\frac{A}{\Jac(A)},
\frac{A}{ \Jac(A)})=0$.
The following extends \cite[Remark 4.6]{A}.

\begin{corollary}
Let $A$ be a  domain of prime characteristic. Let $\frak a$ and $\frak b$ be  radical ideals  of $A^\textbf{+}$. Then $\Tor^{A^\textbf{+}}_i (A^\textbf{+}/
\frak \fa, A^\textbf{+}/ \frak \fb)=0$ for all $i\geq1$.
\end{corollary}

\begin{proof} This is combination of Observation \ref{ob3} (xiii) with the  miraculous vanishing formula.\end{proof}

This is not true for any ideals:

\begin{example}
Let $A$ be a noetherian henselian local domain and let $I$ be a finitely generated proper and nonzero ideal of $A^\textbf{+}$. Then $\Tor^{A^\textbf{+}}_1 (A^\textbf{+}/
I, A^\textbf{+}/ I)\neq0$.
\end{example}

\begin{proof} Recall that  $\Tor^{A^\textbf{+}}_1 (A^\textbf{+}/
I, A^\textbf{+}/ I)\simeq I/I^2$. We need to show $I\neq I^2$. In view of Lemma \ref{art},
 $A^\textbf{+}$ is quasilocal. The claim $I\neq I^2$  follows by Nakayama's lemma.
\end{proof}

However, there is a weak version of rigidity of $\Tor$:

\begin{example}
Let $R$ be a noetherian complete local regular ring of prime characteristic and let $I$  and $J$ be ideals of $R$. If $\Tor_i^{R^\textbf{+}}(\frac{R^\textbf{+}}{IR^\textbf{+}},\frac{R^\textbf{+}}{JR^\textbf{+}})=0$, then $\Tor_j^{R^\textbf{+}}(\frac{R^\textbf{+}}{IR^\textbf{+}},\frac{R^\textbf{+}}{JR^\textbf{+}})=0$  for all $j>i$.
\end{example}

\begin{proof}
By \cite[6.7]{HH1}, $R^\textbf{+}$ is  flat over $R$. Thus $\Tor_i^{R^\textbf{+}}(\frac{R^\textbf{+}}{IR^\textbf{+}},\frac{R^\textbf{+}}{JR^\textbf{+}})\simeq\Tor_i^{R}(R/I,R/J)\otimes_RR^\textbf{+}$. The extension $R\to R^\textbf{+}$ is faithful. Hence, $\Tor_i^{R}(R/I,R/J)=0$.
By the rigidity of $\Tor$-modules over $R$, $\Tor_j^{R}(R/I,R/J)=0$ for all $j>i$ (see \cite{au}). So, $\Tor_j^{R^\textbf{+}}(\frac{R^\textbf{+}}{IR^\textbf{+}},\frac{R^\textbf{+}}{JR^\textbf{+}})\simeq\Tor_j^{R}(R/I,R/J)\otimes_RR^\textbf{+}=0$  for all $j>i$.
\end{proof}

One may prove the following by the straightforward arguments:

\begin{corollary}
Let $R$ be a noetherian regular local ring   of prime characteristic. Let $A$ be any perfect algebra over $R^\infty$. If $A$ is finitely
presented over $R^\infty$, then $A$  is free.
\end{corollary}

\begin{proof}
 In the light of the miraculous vanishing formula we see that $\Tor^{R^\infty}_{>0} (R^\infty/\fm_{R^\infty}, A)=0.$\footnote{Vanishing of the first tor is enough (see \cite[Corollary 2.5.10]{G}). Despite of this, we would like to follow this proof.}The claim  follows from the following fact: \begin{enumerate}
\item[Fact A]: (See \cite[Theorem 3.1.2]{G}) Let $B$ be a coherent ring and let $J\subset \Jac(B)$. If $M$
is a finitely presented $B$-module such that $\Tor_{>0}^B(B/J,M) = 0$, then $\pd_B(M) = \pd_{B/J}(M/JM)$.
 \end{enumerate}
\end{proof}

Let $B\twoheadleftarrow A\twoheadrightarrow C $ be a surjective diagram of perfect rings of prime characteristic.  Let us compute $\Ext_{A}^{>0} (B,C)$ via an example.

\begin{example}
Let $R$ be a 1-dimensional complete local domain of prime characteristic and let $A:=R^\textbf{+}$.  Look at the diagram $R^\textbf{+}/\fm_{R^\textbf{+}}=:B\twoheadleftarrow A\twoheadrightarrow C:= R^\textbf{+}/\fm_{R^\textbf{+}}$ of perfect rings. Then
$\Ext_{A}^i (B,C)=0$ for all $i>0$.
\end{example}

\begin{proof}It is not difficult to see that $\pd(R^\textbf{+}/\fm_{R^\textbf{+}})=2>1=\id(R^\textbf{+}/\fm_{R^\textbf{+}}),$ see \cite[Example 8.2(ii)]{A}). So, the only challenging $\Ext^i$ is $\Ext^1$.
We will use the fact that $(A,\fm_A)$ is a valuation ring. In particular, any finitely generated ideal of $A$ is principal.
Let $I=xA$ be any finitely generated ideal of $A$. Then $$0\lo A \stackrel{x}\lo A\lo A/I\lo 0$$ is a free resolution of $A/I$. By applying
$\Hom(-,A/ \fm_A)$ to it we get to the exact sequence $$ \Hom_A(A,A/ \fm_A) \stackrel{x}\lo \Hom_A(A,A/ \fm_A)\lo\Ext^1_A(A/I,A/ \fm_A) \lo \Ext^1_A(A,A/ \fm_A)=0.$$Note that $\coker (A/ \fm_A \stackrel{x}\lo A/ \fm_A)=A/ \fm_A$ provided $x\in\fm_A$. Thus, $$\Ext^1_A(A/I,A/ \fm_A)=\coker\left(\Hom_A(A,A/ \fm_A) \stackrel{x}\lo \Hom_A(A,A/ \fm_A)\right)=A/ \fm_A \quad(\ast)$$
Let $I=xA$ and $J=yA$ be finitely generated ideals of $A$. Without loss of the generality we may assume that  $t=y/x\in A$, since $A$ is a valuation domain. The natural map $A/I\to A/J\to 0$ induces the map $\Ext^1_A(A/J,A/ \fm_A)\to \Ext^1_A(A/I,A/ \fm_A)$,  a multiplication by $t$.  If $I\neq J$ then $t$ is not invertible. So, $t\in \fm_A$. It turns out that the maps in the inverse system $$A/ \fm_A\longleftarrow A/ \fm_A\longleftarrow \ldots$$ are the zero maps.
Note that $A/\fm_A\simeq\varinjlim A/ I_i$ where $I_i$ is finitely generated. By \cite[Corollary 2.4]{lim}
$\Ext^1_A(\varinjlim M_i,A/ \fm_A)\simeq\vpl\Ext^1_A (M_i, A/ \fm_A)$.
We apply this along with  $(\ast)$ to observe that
$$\Ext^1_A (A/\fm_A, A/\fm_A)\simeq\Ext^1_A (\varinjlim A/ I_i, A/\fm_A)\simeq \vpl\Ext^1_A ( A/ I_i, A/\fm_A)\simeq\vpl A/\fm_A\simeq 0.$$
\end{proof}

\begin{example}
Let $R$ be a 1-dimensional complete local domain of residue prime characteristic and let $A:=R^{\infty}$.  Look at the diagram $R^{\infty}/ \fm_{R^\infty}=:B\twoheadleftarrow A\twoheadrightarrow C:= R^{\infty}$ of perfect rings. Then  $\Ext_{A}^1 (B,C)=0$.
\end{example}

\begin{proof} Let $x\in R$ be the uniformazing element.  Recall that  a free resolution of $B$
over $A$ is given by the following:$$\begin{CD}
0 @>>> \bigoplus _\mathbb{N}R^{\infty} @>X>> \bigoplus _\mathbb{N}R^{\infty} @>Y>>  R^{\infty} @>>>R^{\infty}/\fm_{R^\infty}@>>>0 \quad(\ast)
\\
\end{CD}
$$
where the matrixes $X$ and $Y$ are defined by:
\begin{equation*}
X := \left(
\begin{array}{cccccc}
1      & 0   & 0 & 0 & 0 &\cdots\\
- x^{\frac{p-1}{p}}& 1   & 0 & 0 & 0 &\cdots \\
0 & -x^{\frac{p-1}{p^{2}}}   & 1 & 0 & 0 &\cdots \\
0 & 0   & -x^{\frac{p-1}{p^{3}}}  & 1 & 0 &\cdots\\
0 & 0   & 0 &  -x^{\frac{p-1}{p^{4}}} & 1 & \cdots\\
\vdots & \vdots   &  \vdots& \vdots& \vdots & \ddots\\
\end{array}  \right)
\end{equation*}
and
\begin{equation*}
Y:= \left(
\begin{array}{ccccc}
x& x^{1/p}& x^{1/p^2} & \ldots  \\
\end{array} \right)^t.
\end{equation*}

Note that $$\Hom_{R^{\infty}}(\bigoplus _\mathbb{N}R^{\infty},R^{\infty})\simeq\prod _\mathbb{N}\Hom_{R^{\infty}}(R^{\infty},R^{\infty})\simeq\prod _\mathbb{N} R^{\infty}$$
Anther this identification and after applying $\Hom_{R^{\infty}}(-, R^{\infty})$ to $(\ast)$ we have
$$\begin{CD}
0 @>>> R^{\infty} @>\alpha>>\prod _\mathbb{N} R^{\infty}  @>\beta>> \prod _\mathbb{N} R^{\infty} @>>>0\quad(\ast,\ast)
\\
\end{CD}
$$
where $\alpha$ is the assignment via $a\mapsto (ax^{1/p^n})_{n\in \mathbb{N}_0}$ and $\beta$   assigns the sequence $(a_n)_{n\in \mathbb{N}}$ to the $$\beta(a_n):=(a_1-x^\frac{p-1}{p}a_2,\ldots, a_n-x^{\frac{p-1}{p^{n}}}a_{n+1},\ldots).$$
The homology of $ (\ast,\ast)$ in middle is $\Ext_{A}^1 (B,C)$.
Let $(a_n)\in\ker(\beta)$. Then $0=\beta(a_n)$.
One may read as  $a_{n}=x^{\frac{p-1}{p^{n}}}a_{n+1}$ for all $n$. Define $a(n):=\sum_{i=1}^n\frac{p-1}{p^{i}}$.  Iterate this inductively,  we have $a_1=x^{a(n)}a_n$.
There is a valuation map $v:A\to \mathbb{Q}$ such that if $v(a)\geq1$ then $x\mid a$. Indeed, $R:= V$ is a $\DVR$. Let $v$ be a value map of $V$. Let $r\in R^\infty$. Then $r^{p^n}\in R$ for some $n$. The assignment $r\mapsto v(r^{p^n})/p^n$ defines a normalized value map on $R^\infty$. Since $\lim_{n\to\infty} a(n)=1$ we observe that $v(a_1)\geq1 $. So,
$x\mid a_1$. Set $a:=\frac{a_1}{x}\in A$. Then $\alpha(a)=(a_n)$.  From this we observe $\im(\alpha)=\ker(\beta)$. Therefore, $\Ext_{A}^1 (B,C)=0$.

\end{proof}

We close this section by  extending a result of Bhatt and Scholze  \cite[Proposition 11.29]{BS}.
To this end, we recall the following trick of Auslander \cite[Proposition 3]{au2}.

\begin{lemma}\label{Aus} Let $A$ be a ring and
let $\Gamma$ be a well-ordered set. Suppose that $\{N_\gamma :
\gamma\in \Gamma\}$ is a collection of submodules of an $A$-module
$M$ such that $\gamma'\leq \gamma$ implies $N_{\gamma'}\subseteq
N_{\gamma}$ and $M = \bigcup_{\gamma\in \Gamma} N_\gamma$. Suppose
that $\pd_A (N_\gamma/\bigcup_{\gamma'<\gamma} N_{\gamma'})\leq n$
for all $\gamma\in \Gamma$. Then $\pd_A(M)\leq n$.
\end{lemma}

 \begin{proposition}
 Let $R \to S$ be a perfectly finitely presented map of perfect $\mathbb{F}_p$-algebras. Then $S$ is of
finite  projective dimension over $R$.
\end{proposition}

\begin{proof}
We will show that if $S$ is the perfection of $\frac{R[X_1,\ldots,X_n]}{(f_1,\ldots,f_m)}$, then  $\pd_R(S)\leq2m$. Set $T:=(R[\underline{X}])^\infty$ and let $M$ be any $R$-module. Also, set $(\underline{f}^\infty):=\sum_{i=1}^m (f_i)^\infty$. We look at the base change spectral sequence
$$\Ext^{i}_{T}(T/(\underline{f}^\infty),\Ext_{R}^j(T,M))\Rightarrow\Ext^{i+j}_{R}(T/(\underline{f}^\infty),M).$$Since $T$ is free over $R$ the spectral sequence
collapses and so$$\Ext^{n}_{T}(T/(\underline{f}^\infty),\Hom_{R}(T,M))\simeq\Ext^{n}_{R}(T/(\underline{f}^\infty),M).$$
If we show $\Ext^{i}_{T}(T/(\underline{f}^\infty),-)=0$ for all $i>2m$, then we get the claim. After replacing $R$ with $T$ we may assume that $n=0$.
We do induction by $m$. First, we assume that $m=1$. By Lemma \ref{Aus} $\pd_R(R/(f_1^\infty))\leq2.$ Let $\textbf{Q}^i$ be the corresponding free resolution of length two. By using induction on $m$, we will show that
$\bigotimes_{i=1}^{m} \textbf{Q}^i$  is a  free
resolution of $\frac{R}{(\underline{f}^\infty)}$  of length  at most $2m$. Set $(\underline{g}^\infty):=\sum_{i=2}^{m}(f_i^{\infty})$
By the induction hypothesis, $\textbf{P}:=\bigotimes_{i=2}^{m}
\textbf{Q}^i$ is a free resolution of
$\frac{R}{(\underline{g}^\infty)}$ of length  at most $2m-2$. Recall that
$$H_n(\textbf{Q}^1\otimes_{R}\textbf{P})\simeq \Tor^{R}_n(R/
(f_1^{\infty}),\frac{R}{(\underline{g}^\infty)}).$$ The ideal
$(f_1^{\infty})$ is  flat. Hence, for each $n>1$ we have
$\Tor^{R}_n(R/ (f_1^{\infty}),\frac{R}{
(\underline{g}^\infty)})=0$. Also,
\begin{enumerate}
\item[i)] $\Tor^{R}_1(R/ (f_1^{\infty}),\frac{R}{(\underline{g}^\infty)})
\simeq \frac{(f_1^{\infty})\cap(\underline{g}^\infty)}{
(f_1^{\infty})(\underline{g}^\infty)}\stackrel{\ref{ob3}}=0$,
\item[ii)] $H_0(\textbf{Q}^1\otimes_{R}\textbf{P})\simeq
R/(f_1^{\infty})\otimes_{R}\frac{R}{(\underline{g}^\infty)}\simeq
R/ (\underline{f}^\infty).$
 \end{enumerate}
 Thus,
$\bigotimes_{i=1}^{m}
\textbf{Q}^i$ presents a
free resolution of $R/(\underline{f}^\infty)$ of length   $2m$.
\end{proof}

\section{A glimpse  through mixed-characteristic rings}
In his seminal paper Witt proved that the category of perfect quasilocal rings of prime characteristic $p$ is the same
as of the category of $p$ torsion-free, $p$-adically complete semiperfect quasilocal rings of mixed characteristic $p$.  This result extended by many mathematician. Please see
 Scholze's thesis. Shimomoto asked:

\begin{question}(See \cite[Question 2]{sss})
Let $A$ be a coherent perfect $\mathbb{F}_p$-algebra. Is $\W(A)$  coherent?
\end{question}

And he said: ``This question seems a bit subtle,
because it is known that a power series ring over a valuation ring of Krull dimension greater than
one is not coherent.'' In fact, by using a $p$-adic modification of the quotation, we answer the question.

\begin{observation}\label{obval}
 Let $A$ be a perfect domain of prime characteristic. Then its fraction field $F$ is a perfect field.
Conversely, let $F$ be a perfect field of characteristic $p$ equipped with a valuation $v$ and value ring $A$.
Then $A$ is a perfect ring.
\end{observation}

\begin{proof}
For the first assertion, let $x/y\in F$ be nonzero. The elements  $x^{1/p^n}$  and $y^{1/p^n}$ are in $A$.
So, $x^{1/p^n}/y^{1/p^n}\in F$. Its $p^n$-power is $x/y$. Thus $F$ is perfect. Conversely, let $F$ be a perfect field of characteristic $p$ equipped with a valuation $v$ and value ring $A$. Note that $A$ consists of elements of nonzero value.
Let $a\in A$ be nonzero. Then $a^{1/p^n}\in F$, because $F$ is perfect. Since $v(a^{1/p^n})=1/p^nv(a)\geq0$,
we get that $a^{1/p^n}\in A$ for all $n>0$. This shows that $A$ is perfect.
\end{proof}

\begin{example}
Let $R$ be any $1$-dimensional  local ring of prime characteristic. Then $\W(R^\infty)$ is not coherent
without any regard with respect to coherent property of $R^\infty$.
\end{example}

\begin{proof}Suppose first that $R^\infty$ is  coherent.
By applying Corollary \ref{1}, $R^\infty$ is a valuation ring. This is perfect and coherent.
By Observation \ref{obval}  $Q(R^\infty)$ is a perfect field. Its value group is the set of rational numbers
whose denominator is a $p$-power. In particular, the value group is an strict subgroup of $\mathbb{R}$. Let us follow Kedlaya:
In the light of  \cite[Theorem 1.2]{ko} we see that $\W(R^\infty)$ is not coherent.

Finally, suppose that $R^\infty$ is not coherent.
Suppose on the contradiction that $\W(R^\infty)$ is  coherent.
There is a natural isomorphism $\W(R^\infty)/p\W(R^\infty)\simeq R^\infty$.
In view of \cite[Theorem 2.4.1(1)]{G}, quotient of a coherent ring with a finitely generated ideal
is again coherent. From this we observe that $R^\infty$ is  coherent. This contradiction implies that
$\W(R^\infty)$ is not coherent.
\end{proof}

\begin{fact}(Auslander-�Buchsbaum, P. Kohn, Vasconcelos; see \cite[Theorem 5.1]{V2})\label{v8}
Let $A\to B$ be a ring homomorphism such that $d:=\pd_A(B)<\infty$ and there is an exact complex
$0\to P_d\to\ldots\to P_0\to B\to 0$ of finitely generated projective $A$-modules such that $\Ext^i_A(B,A)=0$ for all
$i<d$. Then for any $B$-module $M$ with $\pd_B(M) < \infty$ we have
$\pd_A(M)=\pd_B(M)+\pd_A(B)$.
\end{fact}

The following result  extends \cite[Lemma 7.8]{BS} by Bhatt and Scholze:

\begin{corollary}
 Let $R$ be a perfect $\mathbb{F}_p$-algebra, and let $Q$ be (not necessarily finitely generated and not necessarily projective) $R$-module of finite projective dimension.
Then $\pd_{\W(R)} (Q)=\pd_{R} (Q)+1.$
\end{corollary}

\begin{proof}
There is an exact sequence $0\to\W(R)\stackrel{p}\lo\W(R)\to\W(R)/p\W(R)\to 0$.
This implies that $\pd_{\W(R)}(\frac{\W(R)}{p\W(R)})=1$ and that $\Hom_{\W(R)}(\frac{\W(R)}{p\W(R)},\W(R))=0$.
In particular, we are in the situation of Fact \ref{v8}. By applying Fact \ref{v8}
we see $\pd_{\W(R)} (Q)=\pd_{R} (Q)+1,$ as claimed.
\end{proof}

Let $A$ be a commutative ring and let $p$ be a non-unit prime  in $A$.
By the \textit{Fontaine ring} of $A$, we mean $$\E(A):=\vpl\xymatrix{(\ldots \ar[r]^{\F}\ar[r]^{\F}&A/pA \ar[r]^{\F}&A/pA)}.$$ If $\nil(A)^n=0$ for some $n\in \mathbb{N}$, then $\E(A)=\E(A_{\red})$. Also, $\E(A)$ is perfect: the $p^{th}$
root of $(r_n)$ is $(s_n)$, where $s_n := r_{n+1}$.

\begin{observation}\label{comea}Here, we  show
$\Wdim(\E(A))$ has properties both similar to and different from those of $\Wdim(A^\infty)$ via some examples.
\begin{enumerate}
\item[i)] (Witt) We look at $A:=\mathbb{Z}_p$ the ring of $p$-adic integers. By Fermat's little theorem, the  Frobenius map is identity over $\mathbb{F}_p$. Then
 $\E(A)=\vpl\xymatrix{(\ldots \ar[r]^{=}\ar[r]^{=}&\mathbb{F}_p \ar[r]^{=}&\mathbb{F}_p)}\simeq\mathbb{F}_p.$ Thus, $$\Wdim(\E(A))=0<1=\Wdim(A).$$

\item[ii)] (This extends \cite[Lemma 3.4(iv)]{Sch} by the same proof) If $A$ is a perfect reduced ring of  characteristic $p$, then  $\E(A)=\vpl\xymatrix{(\ldots \ar[r]^{\simeq}\ar[r]^{\simeq}&A \ar[r]^{\simeq}&A)}\simeq A.$ Thus, $$\Wdim(\E(A))=\Wdim(A).$$

\item[iii)] Let us give an example such that $\Wdim(E(A))=0<n=\Wdim(A),$ where $n\in \mathbb{N}\cup\{\infty\}$. To this end, let $(A,\fm,k)$ be a noetherian complete local ring (not necessarily regular) with perfect residue
field of characteristic $p$. Then$$\E(A)=\E(\vpl A/\fm^i)\simeq \vpl \E( A/\fm^i)\simeq  \vpl \E\left((A/\fm^i)_{ \red}\right)\simeq\vpl \E( A/\fm)\stackrel{(ii)}\simeq \vpl k\simeq k.$$Note that weak dimension of $A$ can be any thing. From this we get the claim.

\item[iv)] Let us give an example such that $\Wdim(E(A))=1<\infty=\Wdim(A)$. To this end, let  $A:=\frac{\mathbb{F}_p[[X]]^\infty}{(X)}$.
By Corollary \ref{1}, $\mathbb{F}_p[[X]]^\infty$ is coherent. In view of \cite[Theorem 2.4.1(1)]{G}, $A$ is coherent. If $A$ were be of finite weak dimension it should be reduced. But $A$ is not reduced, because $(x^{1/p})^p=x=0$. Thus $\Wdim(A)=\infty$. The following completion is the $(X)$-adic. By definition, $$\E(A)=\vpl_{\F} (\frac{\mathbb{F}_p[[X]]^\infty}{(X)})\simeq\vpl _{n\in\mathbb{N}} \left(\frac{\mathbb{F}_p[[X]]^\infty}{(X^{p^{n+1}})}\twoheadrightarrow\frac{\mathbb{F}_p[[X]]^\infty}{(X^{p^n})}\right)\simeq(\mathbb{F}_p[[X]]^\infty)\widehat{ }. $$  Note that such a  completion of a valuation domain is  again valuation domain. Thus, $$\Wdim(E(A))=1<\infty=\Wdim(A).$$

\item[v)] By a result of Gabber and Ramero \cite{GAB}, if $A$ is a valuation domain of mixed  characteristic then $\E(A)$ is a valuation domain. In particular, $\Wdim(\E(A))\leq\Wdim(A)=1$.
\end{enumerate}
\end{observation}

When is $\Wdim(\E(A))<\infty$?
Of course, this is not true in  general:

\begin{example}
Let  $(R,\fm_R)$ be the localization of $\mathbb{F}_p[X_1,\ldots]$ at $(X_1,\ldots)$.  Let $A:=R^\infty$.  In view of Example \ref{needno}
$\Wdim(A)=\infty$. By Observation \ref{comea}(ii), $\E(A)\simeq A$. Thus, $\Wdim(\E(A))=\infty$. Also, $\Wdim(\W(\E(A)))=\infty$.
\end{example}

 \begin{question}(Shimomoto)
How can determine $\Wdim(\E(A))$ in terms of $\Wdim(A)$?
\end{question}

Let $A$ be a mixed characteristic valuation domain. Recall that $\E(A)$ is a valuation ring. So, its weak dimension is one. In view of the following formula$$\W(\E(A))/p\W(\E(A))\simeq \E(A),$$
and by applying Fact \ref{v8}, $\Wdim(\W(\E(A)))\geq2$.
Let $[-]:\E(A)\to \W(\E(A))$ be the Teichm\"{u}ller mapping.  Let $x\in \E(A)$ be such that its radical is the maximal ideal. By the natural isomorphism we get that
$p,[x]$ is a regular sequence on  $\W(R^\infty)$ and that $\rad(p,[x])$ is the maximal ideal of $\W(\E(A))$.
From
this and  Fact \ref{cp} we get that $\Kdepth(\W(\E(A)))=2$. We have no data about of its Krull dimension (resp. its prime spectrum).

\begin{acknowledgement}
I would like to thank  everyone who help me running this project. I thank Shimomoto for a number of valuable comments and  encouragement.  Also, Epstein shared his  interesting example with us.
\end{acknowledgement}

%%%%%%%%%%%%%%%%%%%%%%%%%%%%%%%%%%%%%%%%%


\begin{thebibliography}{99}

\bibitem{ab}
I.M. Aberbach and
M. Hochster, \emph{
Finite tor dimension and failure of coherence in absolute integral
closures}, J.  Pure  Appl. Algebra,  {\bf122}, (1997), 171--184.



\bibitem{an}
 Y. Andr\'{e}, \emph{La conjecture du facteur direct}, arXiv:1609.00345.



\bibitem{Ar}M. Artin,  \emph{
On the joins of Hensel rings}, Adv. Math., {\bf7}, (1971), 282�-296.


\bibitem{A}
M. Asgharzadeh, \emph{Homological properties of the perfect and absolute integral closure of Noetherian domains}, Math. Annalen {\bf 348} (2010), 237�-263.


\bibitem{AT}
M. Asgharzadeh and M. Tousi, \emph{On the notion of Cohen-Macaulayness for non-Noetherian rings}, J. Algebra, {\bf {322}} (2009), 2297--2320.


\bibitem{desing}
M. Asgharzadeh, \emph{Desingularization of regular algebras},   J.  Pure  Appl. Algebra, to appear.

\bibitem{AS}
M. Asgharzadeh, and K. Shimomoto,
\emph{Almost Cohen-Macaulay and almost regular algebras via almost flat extensions}, J. Commut. Algebra,
{\bf 4} (2012), no. 4, 445--478.


\bibitem{au2}
M. Auslander, \emph{On the dimension of modules and algebras III}, Nagoya Math. J.  {\bf9}
(1955) 67--77.

\bibitem{au}M.
Auslander, \emph{Modules over unramified regular local rings}, Ill. J. Math.  {\bf {5}}, 631-�647 (1961).

\bibitem{B}{S.F. Barger}, \emph{A theory of grade for commutative rings},
Proc. AMS., {\bf36}, (1972), 365�-368.

\bibitem{Ber} J. Bertin, \emph{
Anneaux $coh\acute{e}rents\ \ r\acute{e}guliers$}, C. R. Acad. Sci.
Paris, $S\acute{e}r$ A-B, {\bf {273}}, (1971).


\bibitem{BS}
B. Bhatt, and P. Scholze,  \emph{Projectivity of the Witt vector affine Grassmannian}, Invent. Math. {\bf{209}} (2017),  329�-423.


\bibitem{Br}H.
Brenner, \emph{Tight closure and plus closure in dimension two}, Amer. J. Math. {\bf128} (2006), no. 2, 531-–539.

\bibitem{BH}
W. Bruns and J. Herzog,  \emph{Cohen-Macaulay rings}, Cambridge University Press {\bf{39}}, Cambridge, (1998).



\bibitem{CE} H. Cartan, and S. Eilenberg, \emph{Homological
algebra}, Princeton University Press,  1956.



\bibitem{GAB}
O. Gabber and L. Ramero, \emph{Foundations of p-adic Hodge theory}, arxiv:math/0409584.

\bibitem{almost}
O. Gabber and L. Ramero, \emph{Almost ring theory},   Springer-Verlag, LNM {\bf{1800}},
Berlin, (2003).

\bibitem{G}
S. Glaz, \emph{Commutative coherent rings}, Springer-Verlag, LNM  {\bf{1371}}, Berlin,  (1989).

\bibitem{sga2}
A. Grothendieck,
\emph{Cohomologie locale des faisceaux coh\'{e}rents et th\'{e}or�mes de Lefschetz locaux et globaux} (SGA 2).
S\'{e}minaire de G\'{e}om\'{e}trie Alg\'{e}brique du Bois Marie  1962. Amsterdam: North Holland Pub. Co. (1968).


\bibitem{g}
M.J. Greenberg, \emph{Perfect closures of rings and schemes}, Proc. Amer. Math. Soc. \textbf{16} (2) (1965) 313–-317.

\bibitem{lim}
 J.L.Hein, \emph{The convertibility of $\Ext^n_R(-,A)$}, Trans. Amer. Math. Soc.  {\bf{195}} (1974), 243--264.

\bibitem{H1}
M. Hochster, \emph{Canonical elements in local cohomology modules and the direct
summand conjecture}, J. Algebra, {\bf84}, (1983), 503--553.



\bibitem{HH1}
M. Hochster,  and C. Huneke, \emph{Infinite integral
extensions and big Cohen-Macaulay algebras}, Ann. of Math.,
{\bf135}(2), (1992), 53--89.


\bibitem{ko}K.S. Kedlaya, \emph{
Some ring-theoretic properties of $A_{inf}$},	arXiv:1602.09016 [math.NT]

\bibitem{ku} E. Kunz, \emph{Characterizations of regular local rings of characteristic p}, Amer. J. Math. {\bf91} (1969), 772--784.

\bibitem{Mat}
H. Matsumura, \emph{Commutative ring theory}, Cambridge Studies in Advanced Math, \textbf{8}, (1986).



\bibitem{No}
D.G. Northcott, \emph{Finite free resolutions}, Cambridge Tracts Math., vol. \textbf{71}, Cambridge Univ. Press, Cambridge,
(1976).



\bibitem{N}D.G. Northcott,  \emph{On the homology theory of general commutative
rings}, J. London Math. Soc., {\bf 36}, (1961), 231--240.

\bibitem{Ol}J.-P.
Olivier,
\emph{Going up along absolutely flat morphisms},
J. Pure Appl. Algebra {\bf30} (1983),  47-59.

\bibitem{O1}
B.L. Osofsky, \emph{Homological dimensions of modules},
CBMS, {\bf12}, 1971.

\bibitem{O2}
B.L. Osofsky, \emph{
A commutative local ring with finite global dimension and zero divisors},
Trans. Amer. Math. Soc. {\bf141} (1969) 377�-385.

\bibitem{PS2}
C. Peskine and L. Szpiro, \emph{Dimension projective finie et cohomologie locale},
Publ. Math. IHES.  {\bf42} (1973),  47-–119.

\bibitem{Sch}
P. Scholze, \emph{Perfectoid spaces}, Publ. Math. IHES. {\bf{116}} (2012), 245--313.

\bibitem{sib}G. Seibert, \emph{Complexes with homology of finite length and Frobenius functors}, J.
Algebra {\bf125} (1989), 278-287.

\bibitem{Ser}
J.-P. Serre, \emph{Groupes proalgebriques}, Publ. Math. IHES. {\bf{7}} (1960) 5-–67.

\bibitem{ss}
K. Shimomoto, \emph{F-coherent rings with applications to tight closure theory},   J.  Algebra, {\bf{338}},  (2011),  24--34.

\bibitem{sss}
K. Shimomoto, \emph{On the Witt vectors of perfect rings in positive characteristic}, Comm. Algebra {\bf{343}} (2015), no. 12, 5328--5342.



\bibitem{wa}
 K.-i. Watanabe and K.-i. Yoshida, \emph{Hilbert-Kunz multiplicity and an inequality between
multiplicity and colength}, J. Algebra {\bf {230}}, (2000), 295--317.


\bibitem{V2}
W.V. Vasconcelos,  \emph{The rings of dimension two}, Lecture Notes in Pure and Applied Mathematics, {\bf{22}} Marcel Dekker, Inc., New York-Basel, (1976).

\end{thebibliography}
\end{document}